\documentclass[11pt]{article}
\usepackage{algorithm}
\usepackage{algpseudocode}
\usepackage{amsmath}
\usepackage{amssymb}
\usepackage{bm}
\usepackage{booktabs}
\usepackage{color}
\usepackage{cite}
\usepackage{epsfig}
\usepackage{epstopdf}
\usepackage{flafter}
\usepackage{float}
\usepackage{ifthen}
\usepackage{multirow}
\usepackage{amsfonts}
\usepackage{makecell} 
\usepackage{pifont}
\usepackage{subfigure}
\usepackage{times}
\usepackage{threeparttable}
\usepackage{url} 
\usepackage{subfigure,graphicx}
\usepackage{float}
\usepackage{color,multirow}
\usepackage{makecell}
\usepackage{setspace}
\usepackage{cite} 
\usepackage{url} 
\usepackage{ifthen} 
\usepackage[T1]{fontenc}
\usepackage{graphicx}
\usepackage{changepage}
\usepackage[colorlinks,linkcolor=red,anchorcolor=gray,citecolor=green,urlcolor=blue]{hyperref}
\usepackage{graphics}
\usepackage{booktabs}

\usepackage{amsfonts}
\usepackage{color}
\usepackage{amssymb}
\usepackage{amsmath}
\usepackage{indentfirst}

\def\sqr#1#2{{\vcenter{\vbox{\hrule height.#2pt
			\hbox{\vrule width.#2pt height#1pt \kern#1pt \vrule width.#2pt}
			\hrule height.#2pt}}}}
\def\signed #1{{\unskip\nobreak\hfil\penalty50
	\hskip2em\hbox{}\nobreak\hfil#1
	\parfillskip=0pt \finalhyphendemerits=0 \par}}
\def\endpf{\signed {$\sqr69$}}

\def\dbR{{\mathop{\rm l\negthinspace R}}}

\def\3n{\negthinspace \negthinspace \negthinspace }
\def\2n{\negthinspace \negthinspace }
\def\1n{\negthinspace }

\def\dbE{\mathbb{E}}
\def\dbF{\mathbb{F}}

\def\ds{\displaystyle}

\def\dbN{{\mathop{\rm l\negthinspace N}}}
\def\dbP{\mathbb{P}}
\def\dbR{{\mathop{\rm l\negthinspace R}}}
\def\={\buildrel \triangle \over =}

\def\resp{{\it resp. }}
\def\a{\alpha}

\def\g{\gamma}
\def\d{\delta}
\def\e{\varepsilon}

\def\k{\kappa}
\def\l{\lambda}

\def\n{\nabla}
\def\si{\sigma}
\def\t{\times}
\def\f{\varphi}
\def\th{\theta}

\def\i{\infty}
\def\ns{\noalign{\ss} }

\def\G{\Gamma}
\def\D{\Delta}

\def\O{\Omega}

\def\cB{{\cal B}}
\def\cC{{\cal C}}

\def\cE{{\cal E}}
\def\cF{{\cal F}}
\def\cG{{\cal G}}

\def\cI{{\cal I}}
\def\cJ{{\cal J}}
\def\cK{{\cal K}}

\def\cM{{\cal M}}
\def\cN{{\cal N}}

\def\cQ{{\cal Q}}

\def\cl{{\cal l}}
\def\no{\noindent}

\def\ss{\smallskip}
\def\ms{\medskip}
\def\bs{\bigskip}
\def\q{\quad}
\def\qq{\qquad}

\def\esssup{\mathop{\rm esssup}}
\def\max{\mathop{\rm max}}
\def\min{\mathop{\rm min}}
\def\exp{\mathop{\rm exp}}
\def\sup{\mathop{\rm sup}}

\def\pa{\partial}

\def\cd{\cdot}
\def\cds{\cdots}

\def\esssup{\hbox{\rm ess$\,$\rm sup$\,$}}

\def\dist{\hbox{\rm dist$\,$}}

\def\as{\hbox{\rm a.s.{ }}}

\def\cl{\overline}

\def\|{\Big |}
\def\({\Big (}
\def\){\Big )}
\def\[{\Big[}
\def\]{\Big]}
\def\bel{\begin{equation}\label}
\def\ee{\end{equation}}
\def\bt{\begin{theorem}}
\def\bcd{\begin{condition}}
	\def\ecd{\end{condition}}
\def\et{\end{theorem}}
\def\bc{\begin{corollary}}
\def\ec{\end{corollary}}
\def\bde{\begin{definition}}
\def\ede{\end{definition}}
\def\bl{\begin{lemma}}
\def\el{\end{lemma}}
\def\bp{\begin{proposition}}
\def\ep{\end{proposition}}
\def\br{\begin{remark}}
\def\er{\end{remark}}
\def\ba{\begin{array}}
\def\ea{\end{array}}
\def\ed{\end{document}}
\def\ns{\noalign{\ms}}
\def\ds{\displaystyle}

\def\square#1{\vbox{\hrule\hbox{\vrule height#1%
	\kern#1\vrule}\hrule}}
\def\rectangle#1#2{\vbox{\hrule\hbox{\vrule height#1%
	\kern#2\vrule}\hrule}}

\font\tenbb=msbm10 \font\sevenbb=msbm7 \font\fivebb=msbm5

\newfam\bbfam
\scriptscriptfont\bbfam=\fivebb \textfont\bbfam=\tenbb
\scriptfont\bbfam=\sevenbb

\newtheorem{lemma}{Lemma}[section]
\newtheorem{remark}{Remark}[section]

\newtheorem{theorem}{Theorem}[section]
\newtheorem{corollary}{Corollary}[section]

\newtheorem{definition}{Definition}[section]
\newtheorem{proposition}{Proposition}[section]
\newtheorem{condition}{Condition}[section]
\newtheorem{assumption}{Assumption}[section]

\makeatletter

\@addtoreset{equation}{section}
\makeatother

\allowdisplaybreaks

\begin{document}
\title{Stability Estimate  for an Inverse Stochastic Parabolic Problem of Determining Unknown Time-varying Boundary\thanks{This work is partially supported by the NSF of China under grants 11971333, 11931011  and 12025105, and by the Science Development Project of Sichuan University under
	grant 2020SCUNL201. }}

\author{  Zhonghua Liao\thanks{School of Mathematics,  Sichuan University,
	Chengdu 610064,  China. E-mail address: zhonghualiao@yeah.net} \  and\  Qi L\"u\thanks{School
	of Mathematics, Sichuan University, Chengdu
	610064, China. E-mail address: lu@scu.edu.cn.}}

\date{}
\maketitle
\begin{abstract}

Stochastic parabolic equations are widely used to model many random phenomena in natural sciences, such as the temperature distribution in a noisy medium, the dynamics of a chemical reaction in a noisy environment, or the evolution of the density of bacteria population. In many cases, the equation may involve an unknown moving boundary which could represent a change of phase, a reaction front, or an unknown population. In this paper, we focus on an inverse problem where the goal is to determine an unknown moving boundary based on data observed in a specific interior subdomain for the stochastic parabolic equation and prove that the unknown boundary depends logarithmically on the interior measurement.  This allows us, theoretically,  to track  and to monitor the behavior of  unknown boundary from observation in an arbitrary interior domain. The stability estimate  is based on a new Carleman estimate  for  stochastic parabolic equations. As a byproduct, we obtain a quantitative unique continuation property for stochastic parabolic equations.

\end{abstract}

\no{\bf 2020 Mathematics Subject Classification}. 35R30, 60H15

\bs

\no{\bf Key Words}. stochastic parabolic equation, quantitative unique continuation, stability, inverse problems, Carleman estimate.

\section{Introduction}
\par To begin with, we introduce some notations   used in  this paper. Let $ T>0 $ and  $(\O, {\cal F}, \dbF, \dbP)$ with
$\dbF\=\{{\cal F}_t\}_{t \geq 0}$ be a complete
filtered probability space on which a  one dimensional standard Brownian motion
$\{W(t)\}_{t\geq 0}$ is defined. 

Let ${\bf H}$
be a Banach space. 
\par $ \bullet $ By $L^2_{\dbF}(0, T; {\bf H})$
we denote  the  space consisting all ${\bf H}$-valued
$\dbF$-adapted process $X(\cdot)$ such that\linebreak $\dbE
|X(\cdot)|_{L^2(0, T; {\bf H})}^2 < +\infty$.
\par $ \bullet $ By
$L^{\infty}_{\dbF}(0, T; {\bf H})$ we denote the
space of all ${\bf H}$-valued
$\dbF$-adapted bounded processes.
\par $ \bullet $ By
$L^2_{\dbF}(\Omega; C([0, T]; {\bf H}))$ we denote the
space of all ${\bf H}$-valued
$\dbF$-adapted continuous processes $X(\cdot)$
with $\dbE |X(\cdot)|^2_{C([0, T]; {\bf H})} < +\infty$.

\par $ \bullet  $ By $ C_\dbF([0,T]; L^2(\O, {\bf H})) $ we denote the space of all ${\bf H}$-valued $\dbF$-adapted processes satisfying $X(\cd): [0,T]\mapsto L_{\cF_T}^2(\O; {\bf H})$ is continuous.
\par All of the aforementioned spaces are Banach spaces equipped with the canonical norm. 

\par Let  $O \subset  \dbR^n $ $ (n\ge 1) $ be a bounded domain. 
Let  $ \{G(t)\}_{t\in[0,T]} $ be a family of bounded domain with $ C^\infty $ boundaries in $ \dbR^n $ $ (n\ge 1) $ which satisfy the following condition.
\begin{condition}\label{con1}
$G(t)\subset O$ for all $t\in [0,T]$, where   $O\subset  \dbR^n $ is a bounded domain; and $ G(0) $ can be mapped into each domain $ G(t) $   by a $ C^\infty $ diffeomorphism $ \tau(t,\cd ) $ and $ \tau(\cd, \cd) $ depend on $ t $ in $ C^1 $ way, i.e.,
$$
\tau (t, \cd ):  G(0)\mapsto G(t)
$$
is $ C^\i $ diffeomorphism and $  \tau(\cd,\cd)\in C^1([0,T]; C(\cl G(0)))$. 
\end{condition}

Put
$$
G((0,T))\= \bigcup_{t\in(0,T)}(\{t\}\t G(t)).
$$ 
We will also denote similar notations like $ \cd((\cd,\cd )) $ for other  set families and time interval in the rest of this paper. 


Assume  $ I(t) \subset \pa G(t)$ $(t\in (0,T)) $   is unknown, which represents, for example, a change of phase, a reaction front, or another unknown population and the rest part of boundaries $ \G \= \pa G(t)\setminus I(t) $ $(t\in (0,T))$ are not time varying. 

Consider the following stochastic parabolic equation on $ G((0,T)) $:
\begin{equation}\label{CD}
\left\{\ba{ll}
\ds du-\D udt=(a_1\cd \n u+b_1 u)dt+c_1udW(t)& \mbox{in } G((0,T)),\\
\ns\ds u=f& \mbox{on }\G,\\
\ns\ds u=0&\mbox{on }I((0,T)),\\
\ns\ds u(0, \cd)=u_0, &\mbox{in } G(0).
\ea \right.
\end{equation} 
Here $a_1\in
L^{\infty}_{\dbF}(0,T;C^{n}(O;\dbR^n))$,
$b_1\in
L^{\infty}_{\dbF}(0,T;C^{n}(O))$,
$c_1\in L^{\infty}_{\dbF}(0,T;
C^{n}(O))$, $ f$ and $u_0 $ denote boundary input and initial input data, respectively. 

The above equation are widely used to model many random phenomena in natural sciences, such as the temperature distribution in a noisy medium, the dynamics of a chemical reaction in a noisy environment, or the evolution of the density of bacteria population (e.g. \cite{PK, QLXZ, GDJZ}).

By coordinate transformation, we can define the (weak) solutions of stochastic parabolic equation (\ref{CD}) by means of  the following variable coefficients parabolic equations defined in cylinder time-space domain:
\begin{equation}\label{dd1}
\ba{ll}
\ds dv-\sum_{i=1}^n\[\langle D_y^2 v\pa_{x_i}\rho(t,x), \pa_{x_i} \rho(t,x)\rangle+\n v\cd \pa_{x_i}^2\rho(t,x)\]dt+\pa_t\rho(t,x)\cd \n vdt\\
\ns\ds =\(a_1\cd \sum_{i=1}^n \pa_{x_i} \rho(t,x) \pa_{y_i} v+a_2v\)dt+c_1vdW(t), \qq  (t,y)=(0,T)\t G(0),
\ea 
\end{equation}
where $ D^2_y $ denote the Hessian matrix of $ v $ in the space variable $ y $, and  $ \rho(t, \cd)=\tau^{-1}(t, \cd)$, $ (t,y)=\rho(t,x) $, $ u(t,x)= v(t,y)$. The uniform ellipticity for the principal part of (\ref{dd1}) is proved in \cite{JLLAMEZ}.  Then the well-posed results for  equation (\ref{dd1}) can be  followed by general  variable coefficient stochastic parabolic equation. See \cite[Chapter 3]{QLXZ} for example.

In the equation \eqref{CD}, we can choose $f$ and $u_0$ according to our needs. The objective  is   to track and monitor $ I((0,T))\= \cup_{t\in[0,T]} \{t\}\t I(t) $ from the special chosen $f$, $u_0 $ and appropriate measurements on the solution $u$. More specifically, we address the following two questions:
\begin{itemize}
\item {\bf Problem 1: Uniqueness.} Can  we determine the behavior of the unknown moving boundaries \linebreak $ \{ I(t)\}_{t\in[0,T]} $  uniquely by using observations from any interior subset $ (0,T) \t O_0\subset G((0,T)) $ under certain particular choices of $f$ and $u_0$?
\item   {\bf  Problem 2: Stability.} If the answer for the first question is affirmative, how does the behavior of the unknown boundaries depend continuously on the observation data? Additionally, what is the explicit convergence rate under certain particular choices of  $f$ and $u_0$? 
\end{itemize}
%


The above problems are inverse problems for stochastic parabolic equations. 
The inverse problems that determining (moving or invariant) unknown boundaries for deterministic parabolic equations arise  in applied contexts of nondestructive testing of materials for either electric or thermal conductors and has been studied extensively (e.g. \cite{GAEBERSV, SV, BCERSV, JEGHMY, KBLC}  and the rich references therein). Nevertheless, as far as we know, there is no published works addressing that problem for stochastic parabolic equations. The primary reason, in our opinion, is that many sharp methods for solving inverse problems of deterministic PDEs cannot be applied to solve the corresponding stochastic problems. For example,   since the sample path of Brownian motion is almost surely non-differentiable for the time variable, the solutions of stochastic parabolic equations is less regular with respect to time.  Unlike the deterministic parabolic equations, many regularity estimates and extension theorems for deterministic case are not applicable for stochastic parabolic equations.  Nevertheless, due to their indispensable applications, these problems have attracted more and more researchers' attention in recent years. Many important progress has been made and many other types of inverse problems have been considered, such as an inverse source problem, determining history data, etc. We refer to \cite{DCYOSP, FDWD, QL1, LZ,  Wu2020, Wu2022, GY, Yuan}  and rich references therein. 

\par Now, we state our main results in this paper. Denote
\begin{equation}\label{a11}
{\bf G}\=\Big\{G((0,T))\|~G((0,T)) \mbox{ satisfy Condition \ref{con1} and }(0,T)\t \G\subset \pa G((0,T))  \Big\}.
\end{equation}
Let $ G_1((0,T)), G_2((0,T))\subset {\bf G} $ with  different unknown boundaries $ I_1((0,T)) $ and $ I_2((0,T)) $ respectively. Suppose $ I(0) $ is known, i.e. $ G_1(0)=G_2(0) $ and the boundary input $f$ and initial input data $ u_0  $ are the same. The observation subset $O_0$ satisfy $ (0,T)\t  O_0\subset G_1((0,T))\cap G_2((0,T)) $. 

We introduce  the following condition about the boundary data.
\begin{condition}\label{assumption4}
The boundary data $ f(\cd,\cd) $ is nontrivial  function on $ \G\times (0,T] $, i.e.,  there exist a positive constant $ F $  such that $ \Vert f(t, \cd)\Vert_{L^\i(\G)}^2\ge F $ for any $ t\in (0,T] $.
\end{condition} 
\begin{remark}
In the study of Problems 1 and 2, it is usual for us to have the freedom to choose the function $f$. Therefore, Condition \ref{assumption4} can always be satisfied as we can make the appropriate choice for $f$. Here we write it as a condition to emphasize our choice of $f$. 

On the other hand, it is important to note that the uniqueness property may not hold even for the deterministic parabolic case without additional conditions on the boundary input data $f$, as shown by the counterexample in \cite{KBLC}.
\end{remark}

Let  $u_i$  be the weak solution   of  (\ref{CD})  on $ G_i((0,T)) $    (with the diffeomorphism $ \tau_i$ ) satisfying $ u_i\circ \tau_i \in   L_\dbF^2 (\O; C([0,T]; L^2(G_i(0))))\cap L_\dbF^2(0,T; H^1(G_i(0)))  $ ($ i=1,2 $). In what follows, we call such solution $ u_i\in L_\dbF^2 (\O; C([0,T]; L^2(G_i(t))))\cap L_\dbF^2(0,T; H^1(G_i(t)))  $ for the convenience of notations. Our first  result is as follow.
\begin{theorem}\label{uniqueness}
Under  Condition \ref{assumption4}, for any $ t_0\in (0,T] $, if 
\begin{equation}\label{o1}
	\dbE \int_0^{t_0}\int_{O_0} |u_1(t, x)-u_2(t,x)|^2dx dt =0, 
\end{equation}
then  $G_1(t_0)=G_2(t_0)$.
\end{theorem}

\par Next, we study the stability estimate for the inverse stochastic parabolic problem with unknown time-varying boundary.  To this end, let us first introduce some  prior assumptions as follows.

\begin{assumption}\label{assumption2}
The solution $ u_i $ ($ i=1,2 $) of (\ref{CD}) satisfying $ u_i\circ \tau_i \in C_\dbF\left([0,T]; L^2(\O; C^1(G_i(0)))\right) $ and there exist a positive constant $\k_0\ge e$ such that $ 
\esssup_{t\in (0,T)} \dbE \Vert u_i(t,\cd)\Vert_{W^{1,\i}(G_i(t))}\le \kappa_0$ ($ i=1,2 $). 
\end{assumption}
We also simply denote the solution in Assumption \ref{assumption2} by $ u_i\in C_\dbF\left([0,T]; L^2(\O; C^1(G_i(t))) \right) $ hereafter. 
\begin{remark}
According to the classical regularity result for stochastic parabolic equations, Assumption \ref{assumption2} can be satisfied when $f$ and $u_1(0)=u_2(0)$ are sufficiently smooth and satisfy certain compatibility relations (e.g.,\cite[Theorem 3.7]{Flandoli1990}). In   Problem 2, both $f$ and $u_1(0)=u_2(0)$ can be chosen at our discretion. Hence, Assumption \ref{assumption2} can be fulfilled. 
\end{remark}

\par In what follows, we denote $ \cB_r(y)\=\{x\in \dbR^n |~|x-y|<r\} $ and $ \cB_r\=\cB_r(0) $ for $r>0$.

\begin{assumption}\label{assumption1}
There exists a constant $ R_0 >0$, such that  the following conditions hold true:
\\	(i) For any $ x\in \pa G_i(t) $,  we have $ \cB_{R_0}(x-R_0\nu(x))\subset G_i(t) $, where $ \nu(x) $ denotes the unit outer normal on $ x $.
\\	(ii) For any $ x_1\in G_i(t) $ and $ x_2\in \pa G_i(t)\cap \cB_{R_0}(x_1)  $, we have $ (x_2-x_1)\cd \nu(x_2)\ge 0 $.
\end{assumption}
\begin{assumption}\label{assumption3}
There exists a constant $ E\ge 1 $ such that for every $ t_0\in (0,T] $, $ R>0 $ and $ x_0\in \dbR^n $ satisfying if $ \cB_{ER}(x_0)\subset G_i(t_0) $, then $ (\max\{t_0-R^2, 0\}, t_0)\t \cB_R(x_0) \subset G_i((0,t_0))$.  
\end{assumption}
\begin{remark} 
Assumption \ref{assumption3} says  that the unknown time-varying boundary $ I_i((0,T)) $ ($ i=1,2 $) should not  change too fast.
\end{remark}

In the following, we will use the symbol $C$ to represent a generic positive constant  depending on $O, T, \G,  a_1, b_1, c_1$, $O_0, F$, $ R_0$ and $E$.  The value of CC may vary from line to line. Additionally, the other absolute constants defined in the remaining part may also depend on the aforementioned parameters.  

Now  we are in the position to present our stability estimate result.
\begin{theorem}\label{theorem1}
Suppose Condition \ref{assumption4} and Assumptions \ref{assumption2}--\ref{assumption3} hold. If  for a given  $\tilde \e\in (0,1)$,   
\begin{equation}\label{e}
	\sup_{t\in (0,T]}\dbE  \int_{O_0}|u_1-u_2|^2 dx\le \tilde \e,
\end{equation}
then  for any $t_0\in (0,T]$, 
\begin{equation}\label{stability}
	{\bf d}( G_1(t_0),  G_2(t_0) )\le  \exp(e^{\g(t_0)^C})|\ln \tilde \e|^{-\g(t_0)^{-C}},
\end{equation}
where $ \g(t)\=(\ln \k_0)^2+e^{\frac{1}{t^{n/2}}} $ and  ${\bf d}(\cd, \cd)$ is the Hausdorff distance defined as
\begin{equation}\label{20-2}
	{\bf d}( G_1(t_0),  G_2(t_0) )=\max\Big\{\sup_{x\in G_1(t_0)} \dist(x, G_2(t_0)), \sup_{x\in G_2(t_0)} \dist(x, G_1(t_0))\Big\}.
\end{equation}
\end{theorem}

In \cite{SV}, the authors solved the problem of determining unknown boundary portions of a thermic conducting body from Cauchy data for a deterministic parabolic equation by establishing a suitable quantitative unique continuation property for that equation. We have taken some ideas from their work. However, there are very few known results for the quantitative unique continuation of stochastic parabolic equations (for example, \cite{ZLQL, XZ, LY, QL1}), and they all seem difficult to apply to our specific problem. Therefore, we first establish a quantitative unique continuation property in the form of a ``two sphere one cylinder inequality" at the moving boundary for stochastic parabolic equations. Our methods are based on a new Carleman estimate for stochastic parabolic operators with both time and space boundary observation terms (see Theorem \ref{theorem_Carleman}) and a fundamental auxiliary lemma (see Lemma \ref{inital_lemma}). Under suitable regularity assumptions, by utilizing the two sphere one cylinder type inequality, we prove some primary estimates and a small propagation estimate. This leads us to obtain the logarithmic stability estimate for determining the state of the unknown boundary problem.



\par The remaining part of this paper is organized as follows. In Section 2, we demonstrate a Carleman estimate for stochastic parabolic equations. Section 3 is dedicated to establishing the quantitative unique continuation results for stochastic parabolic equations. Finally, in Section 4, we present the proofs of Theorems \ref{uniqueness} and \ref{theorem1}.

\section{Carleman estimate for stochastic parabolic equations}\label{car}

In this section, we establish a Carleman estimate for stochastic parabolic equations.

For $ s\in (0, \frac{1}{e}] $, let
\begin{equation}\label{si}
\si(s)\=s\exp\Big\{-\int_0^{ s}\[1-\exp\(-\frac{1}{(\ln t)^2}\)\]\frac{dt}{t}\Big\},
\end{equation}

Based on elementary calculations, we can derive the following straightforward result. For the sake of brevity, we will omit the proof.
\begin{lemma}\label{lemma_si}
For $ s\in (0,\frac{1}{e}] $,  $ \si $ is the solution to the Cauchy problem 
\begin{equation}\label{si1}
	\frac{d}{ds}\ln \(\frac{\si}{s\si'}\)=\frac{2}{s|\ln s|^3} , \q \si(0)=0, \q \si'(0)=1,
\end{equation}
and  there exist a constant $ C_0 >0$ such that 
\begin{equation}\label{si2}
	\ba{ll}
	\ds s e^{-C_0}\le \si(s)\le s,\qq
	e^{-C_0}\le \si'(s)\le 1.
	\ea 
\end{equation}
\end{lemma}

\par Here is the main result of this section.
\begin{theorem}\label{theorem_Carleman}
Let $ (t_0,x_0)\in G((0,T)) $ with the time-varying boundary $ I((0,T)) $. Let 
\begin{equation}\label{u}
	\ba{ll}
	\ds  u\in L_\dbF^2(\O; C([0,T]; L^2(G (t))))\cap L_\dbF^2(0,T; H_0^1(G (t))),
	\ea 
\end{equation}
satisfies $ u\equiv 0 $ on $ \pa G(t) $ $ (t\in (0,T)) $ and solves
\begin{equation}\label{u1}
	du-\D udt=g_1 dt+g_2dW(t)
\end{equation}
for some $ g_1\in L_\dbF^2(0,T; L^2(G(t))) $ and $ g_2\in L_\dbF^2(0,T; H^1(G(t))) $. For any $\l\ge 1$, $ a\in (0,+\i) $, $b \in (0,t_0]$  and $ a+b\le \frac{1}{e} $,  we have the following inequality
\begin{equation}\label{Carleman}
	\ba{ll}
	\ds \dbE\int_{Q_b}   \frac{\l e^{-C_0}}{|\ln (t_0-t+a)|^3} e^{2\phi}u^2dxdt + \dbE \int_{Q_b} \frac{\frac{1}{2}e^{-2C_0}\si_a}{|\ln (t_0-t+a)|^3}e^{2\phi} |\n u|^2dxdt\\
	\ns\ds \le e^{2C_0}\dbE \int_{Q_b} \si_a^{2}e^{2\phi} g_1^2dxdt+\cM+\cN,
	\ea 
\end{equation}
where $Q_b\= G((t_0-b, t_0))$,  $C_0$ is given in (\ref{si2}),
\begin{equation}\label{a1}
	\ba{ll}
	\ds \si_a(t)\= \si(t_0-t+a), \q \phi(t,x)\=-\frac{|x-x_0|^2}{8(t_0-t+a)}-\l \ln \si_a, 
	\ea 
\end{equation}
\begin{eqnarray}\label{cB}
&&\cM = \dbE \int_{G(t_0-b)} (a+b)^2|\n( e^{\phi} u(t_0-b,x))|^2dx+\dbE\int_{G(t_0)}a \l e^{C_0}\si(a)^{-2\l} u(t_0,x)^2dx \nonumber\\
&& \qq\q +\dbE \int_{G(t_0-b)} \(\frac{|x-x_0|^2}{16}+\frac{a+b}{2}\)\si(a+b)^{-2\l}u(t_0-b,x)^2dx\\
&& \qq\q  -\dbE \int_{t_0-b}^{t_0} \int_{\pa G(t)} (x-x_0)\cd \nu \(|\n v|^2+\frac{1}{4}(t_0-t+a)e^{2\phi}g_2^2\)dxdt, \nonumber
\end{eqnarray}
and
\begin{equation}\label{cB1}
	\ba{ll}
 \ds \cN= \dbE \int_{Q_b } \(\frac{1}{2}+\frac{n}{4}\)(t_0-t+a)e^{2\phi} g_2^2dxdt+\dbE \int_{Q_b} (t_0-t+a)^2e^{2\phi} |\n g_2|^2dxdt. 
	\ea
\end{equation}
\end{theorem}
%


\begin{remark}
The construction of the weight function $\phi$ in equation (\ref{a1}) is based on the idea from \cite{LEFJF}, where the authors aim to prove strong unique continuation for deterministic parabolic equations. However, we need to make adjustments to the form of $\sigma$ in order to ensure the applicability of the Carleman estimate in proving the two-sphere one-cylinder type inequality for stochastic parabolic equations in Section 3.
\end{remark}

\par {\bf Proof of Theorem \ref{theorem_Carleman}.} 
Let $v\=e^{\phi} u$. Then we have
\begin{equation}\label{p1}
e^{\phi} (du-\D u dt)= dv-\phi_t vdt-(\D v+|\n \phi|^2 v-2\n\phi\cd \n v-\D \phi v)dt.
\end{equation} 
Consequently,
\begin{equation}\label{p2}
(t_0-t+a)e^{\phi} (du-\D udt)=\cI_1dt+\cI_2,
\end{equation}
with
\begin{equation}\label{p3}
\begin{cases}
	\ds \cI_1\= -(t_0-t+a)\(\D v+(|\n \phi|^2+\phi_t)v\)+\frac{1}{2}v,\\
	\ns\ds \cI_2\= (t_0-t+a)\[dv+\(\D \phi v+2\n \phi\cd\n v\)dt\]-\frac{1}{2}vdt.
\end{cases}
\end{equation}
Noting that $ \n\phi =-\frac{x-x_0}{4(t_0-t+a)} $ and $ \D \phi=-\frac{n}{4(t_0-t+a)} $, we have
\begin{equation}\label{phi}
\cI_2= (t_0-t+a)\[dv-\( \frac{n}{4(t_0-t+a)}v+\frac{(x-x_0)\cd \n v}{2(t_0-t+a)}\)dt\]-\frac{1}{2}vdt.
\end{equation}
It follows from (\ref{p2}) and (\ref{p3}) 
that
\begin{equation}\label{p4}
\ba{ll}
\ds \dbE \int_{Q_b} \cI_1^2+(t_0-t+a)^2e^{2\phi}g_1^2dxdt\\
\ns\ds \ge \dbE \int_{Q_b}2\cI_1(t_0-t+a)e^\phi (du-\D udt)dx\\
\ns\ds \ge \dbE \int_{Q_b}\[ 2\cI_1^2dxdt +2\cI_1\cI_2\].
\ea 
\end{equation}
Then, by directly calculation, we have
\begin{equation}\label{I1I2_1}
2\cI_1\cI_2=\cJ_1+\cJ_2+\cJ_3,
\end{equation}
where
\begin{equation}\label{I1I2_2}
\begin{cases}
	\ds \cJ_1\= \frac{1}{2}(t_0-t+a) \(\D v+(|\n \phi|^2+\phi_t)v\)\(nv+2(x-x_0)\cd\n  v\)dt,\\
	\ns\ds \cJ_2\= -(t_0-t+a)\(\D v+(|\n \phi|^2+\phi_t)v\)\(2(t_0-t+a)dv-vdt\),\\
	\ns\ds \cJ_3\= (t_0-t+a)vdv-\(\frac{1}{2}+\frac{n}{4}\)v^2dt-\frac{1}{2}v(x-x_0)\cd\n vdt.
\end{cases}
\end{equation}

Next, we estimate $ \cJ_1, \cJ_2$ and $\cJ_3  $, respectively. For $ \cJ_1 $, noting that
\begin{equation}\label{2}
\phi_t= -\frac{|x-x_0|^2}{8(t_0-t+a)^2}-\l \frac{\pa_t \si_a}{\si_a},
\end{equation}
where $ \pa_t\si_a(t)=-\si'(s)|_{s=t_0-t+a} $,   we have
\begin{equation}\label{J1_1}
\cJ_1=\frac{1}{2}(t_0-t+a)\(\D v+\(-\frac{|x-x_0|^2}{16(t_0-t+a)^2}-\l \frac{\pa_t\si_a}{\si_a}\) v \)(nv+2(x-x_0)\cd\n v)dt.
\end{equation}
By direct computations, we have
\begin{equation}\label{J1_2}
\ba{ll}
\ds \frac{1}{2}(t_0-t+a)\D v(n v+2(x-x_0)\cd \n v)\\
\ns\ds =\frac{1}{2}(t_0-t+a) \[n\n\cd(\n v v) -n|\n v|^2+2\n\cd (\n v(x-x_0)\cd \n v)-2|\n v|^2\\
\ns\ds \qq\qq\qq\q -\n\cd  ((x-x_0)|\n v|^2)+n|\n v|^2\]\\
\ns\ds = \frac{1}{2}(t_0-t+a)\[\n\cd \(nv\n v+2\n v(x-x_0)\cd \n v-(x-x_0)|\n v|^2\)-2|\n v|^2\],
\ea 
\end{equation}
and
\begin{eqnarray}\label{J1_3}
&& \frac{1}{2}(t_0-t+a)\(-\frac{|x-x_0|^2}{16(t_0-t+a)^2}-\l \frac{\pa_t \si_a}{\si_a}\) v(nv+2(x-x_0)\cd\n v) \nonumber\\
&& =\n\cd \[\frac{1}{2}(t_0-t+a)(x-x_0) v^2 \(-\frac{|x-x_0|^2}{16(t_0-t+a)^2}-\l \frac{\pa_t\si_a}{\si_a}\)\] \\
&& \q +\(\frac{n}{2}-\frac{n}{2}\)(t_0-t+a)\(-\frac{|x-x_0|^2}{16(t_0-t+a)^2}-\l \frac{\pa_t\si_a}{\si_a}\)v^2 +\frac{|x-x_0|^2}{16(t_0-t+a)}v^2.\nonumber
\end{eqnarray}
Together with (\ref{J1_2}) and (\ref{J1_3}), noting that $ v\equiv 0 $ on $ \pa G(t) $  for any $ t\in (0,T) $, we have
\begin{equation}\label{J1_4}
\ba{ll}
\ds \dbE \int_{Q_b} \cJ_1dx \3n&\ds = \dbE \int_{Q_b} \Big\{-(t_0-t+a)|\n v|^2 +\frac{|x-x_0|^2}{16(t_0-t+a)}v^2\Big\}dxdt\\
\ns&\ds \q +\dbE\int_{t_0-b}^{t_0}\int_{\pa G(t)} (x-x_0)\cd \nu |\n v|^2dxdt.
\ea 
\end{equation}

Next, we hand the term   $ \cJ_2 $. It follows from \eqref{I1I2_2} that
\begin{equation}\label{J2}
\cJ_2=-(t_0-t+a)\[\D v+\(-\frac{|x-x_0|^2}{16(t_0-t+a)^2}-\l \frac{\pa_t \si_a}{\si_a}\) v \]\[2(t_0-t+a)dv-vdt\].
\end{equation}
By It\^o's formula, we find that
\begin{eqnarray}\label{J2_1}
&& -(t_0-t+a)\D v\[2(t_0-t+a)dv-vdt\] \nonumber\\
&& = -(t_0-t+a)\n\cd \[2(t_0-t+a)\n v dv-\n vvdt\]+2(t_0-a+a)^2 \n vd\n v   -(t_0-t+a)|\n v|^2dt \nonumber\\
&& =-(t_0-t+a)\n\cd \(2(t_0-t+a)\n v dv-\n vvdt\)-(t_0-t+a)|\n v|^2dt\\
&& \q +d\((t_0-t+a)^2|\n v|^2\)+2(t_0-t+a)|\n v|^2dt-(t_0-t+a)^2|d\n v|^2 \nonumber\\
&& =(t_0-t+a)\n\cd \(-2(t_0-t+a)\n v dv+\n vvdt\)+d\((t_0-t+a)^2|\n v|^2\) \nonumber\\
&& \q +(t_0-t+a)|\n v|^2dt-(t_0-t+a)^2| \n (e^\phi g_2)|^2dt. \nonumber
\end{eqnarray}
Direct computations give us that
\begin{equation}\label{19-1}
\ba{ll}
\ds -(t_0-t+a)^2|\n (e^{\phi}g_2)|^2dt\\
\ns\ds =-(t_0-t+a)^2e^{2\phi}\(|\n \phi|^2g_2^2+|\n g_2|^2+2g_2\n\phi\cd\n g_2\)dt\\
\ns\ds =-(t_0-t+a)^2\n\cd (e^{2\phi}g_2^2\n \phi )dt+(t_0-t+a)^2e^{2\phi}\(\D \phi g_2^2+|\n \phi|^2 g_2^2-|\n g_2|^2\)dt\\
\ns\ds =-(t_0-t+a)^2\n\cd (e^{2\phi}g_2^2\n \phi )dt -\frac{n}{4}(t_0-t+a)e^{2\phi}g_2^2dt+\frac{|x-x_0|^2}{16}e^{2\phi}g_2^2dt\\
\ns\ds\q -(t_0-t+a)^2e^{2\phi}|\n g_2|^2dt.
\ea 
\end{equation}
This, together with (\ref{J2_1}), implies that 
\begin{equation}\label{19-2}
\ba{ll}
\ds -(t_0-t+a)\D v\(2(t_0-t+a)dv-vdt\)\\
\ns\ds =(t_0-t+a)\n\cd \[-2(t_0-t+a)\n v dv+\n vvdt-(t_0-t+a) \n \phi e^{2\phi}g_2^2dt\]\\
\ns\ds \q +d\[(t_0-t+a)^2|\n v|^2\]+(t_0-t+a)|\n v|^2dt-\frac{n}{4}(t_0-t+a)e^{2\phi}g_2^2dt\\
\ns\ds \q +\frac{|x-x_0|^2}{16}e^{2\phi}g_2^2dt-(t_0-t+a)^2e^{2\phi}|\n g_2|^2dt.
\ea 
\end{equation}
Moreover, we have
\begin{eqnarray}\label{J2_2}
&& -(t_0-t+a)\[-\frac{|x-x_0|^2}{16(t_0-t+a)^2}-\l \frac{\pa_t\si_a}{\si_a}\] v\[2(t_0-t+a)dv-vdt\]\nonumber\\
&&=-d\Big\{(t_0-t+a)^2\[ -\frac{|x-x_0|^2}{16(t_0-t+a)^2}-\l \frac{\pa_t \si_a}{\si_a}\]v^2\Big\}\nonumber\\
&& \q-(t_0-t+a)\(-\frac{|x-x_0|^2}{16(t_0-t+a)^2}-\l \frac{\pa_t\si_a}{\si_a}\) v^2dt \nonumber\\
&& \q +(t_0-t+a)^2\[-\frac{|x-x_0|^2}{8(t_0-t+a)^3}-\l\pa_t\(\frac{\pa_t \si_a}{\si_a}\) \]v^2dt\nonumber\\
&&\q  +(t_0-t+a)^2\[-\frac{|x-x_0|^2}{16(t_0-t+a)^2}-\l \frac{\pa_t \si_a}{\si_a}\]e^{2\phi}g_2^2dt\\
&& =-d\((t_0-t+a)^2\(-\frac{|x-x_0|^2}{16(t_0-t+a)^2}-\l \frac{\pa_t \si_a}{\si_a}\)v^2\)-\frac{|x-x_0|^2}{16(t_0-t+a)}v^2dt\nonumber\\
&& \q +\[-(t_0-t+a)^2\l \pa_t \(\frac{\pa_t\si_a}{\si_a}\)+(t_0-t+a)\l\frac{\pa_t\si_a}{\si_a}\]v^2dt\nonumber\\
&&\q  +(t_0-t+a)^2\(-\frac{|x-x_0|^2}{16(t_0-t+a)^2}-\l \frac{\pa_t\si_a}{\si_a}\)e^{2\phi}g_2^2dt.\nonumber
\end{eqnarray}
Noting $ \pa_t\si_a(t) =-\si'(s)|_{s=t_0-t+a}$, $ \pa_t(\frac{\pa_t\si_a}{\si_a})(t)= (\frac{\si'}{\si})'(s)|_{s=t_0-t+a}$, and  from (\ref{si1}), we obtain that
\begin{equation}\label{th}
-s^2\(\frac{\si'}{\si}\)'-s\frac{\si'}{\si}=s^2\frac{\si'}{\si}\(-\frac{\si''}{\si'}+\frac{\si'}{\si}-\frac{1}{s}\)=s^2\frac{\si'}{\si}\frac{d}{dt}\ln \(\frac{\si}{s\si'}\)=\frac{2s\si'}{|\ln s|^3\si}.
\end{equation}
Combine (\ref{J2_2}) with (\ref{th}), we get
\begin{equation}\label{J2_2*}
\ba{ll}
\ds -(t_0-t+a)\(-\frac{|x-x_0|^2}{16(t_0-t+a)^2}-\l \frac{\pa_t\si_a}{\si_a}\) v\[2(t_0-t+a)dv-vdt\]\\
\ns\ds =d\[(t_0-t+a)^2\(\frac{|x-x_0|^2}{16(t_0-t+a)^2}+\l \frac{\pa_t \si_a}{\si_a}\)v^2\]-\frac{|x-x_0|^2}{16(t_0-t+a)}v^2dt\\
\ns\ds \q -\l\frac{2(t_0-t+a)\pa_t \si_a}{|\ln (t_0-t+a)|^3 \si_a} v^2dt  +\[-\frac{|x-x_0|^2}{16}-\l(t_0-t+a)^2 \frac{\pa_t\si_a}{\si_a}\]e^{2\phi}g_2^2dt.
\ea 
\end{equation}
From (\ref{19-2}) and (\ref{J2_2*}), we find that
\begin{equation}\label{J2_3}
\ba{ll}
\ds \dbE \int_{Q_b}\cJ_2 dx\\
\ns\ds = \dbE \int_{G(t_0)} a^2|\n v(t_0,x)|^2dx-\dbE \int_{G(t_0-b)}(a+b)^2|\n v(t_0-b,x)|^2dx\\
\ns\ds \q +\dbE\int_{G(t_0)}a^2 \(\frac{|x-x_0|^2}{16a^2}+\l\frac{\pa_t\si_a}{\si_a}(t_0)\)v(t_0,x)^2dx\\
\ns\ds   \q -\dbE \int_{G(t_0-b)}(a+b)^2\(\frac{|x-x_0|^2}{16(a+b)^2}+\l\frac{\pa_t\si_a}{\si_a}(t_0-b)\)v(t_0-b,x)^2 dx\\
\ns\ds\q  +\dbE \int_{t_0-b}^{t_0}\int_{\pa G(t)}\frac{1}{4}(t_0-t+a) (x-x_0)\cd\nu e^{2\phi}g_2^2dx \\
\ns\ds \q +\dbE \int_{Q_b}(t_0-t+a)|\n v|^2dxdt -\dbE \int_{Q_b} \frac{|x|^2}{16(t_0-t+a)}v^2dxdt\\
\ns\ds \q -\dbE \int_{Q_b}\l\frac{2(t_0-t+a)\pa_t \si_a}{|\ln (t_0-t+a)|^3 \si_a}  v^2dxdt\\
\ns\ds\q + \dbE \int_{Q_b} (t_0-t+a)^2e^{2\phi}\[-|\n g_2|^2-\l\frac{\pa_t\si_a}{\si_a}g_2^2-\frac{n}{4(t_0-t+a)}g_2^2dt\]dxdt.
\ea 
\end{equation}
\par Combining (\ref{J1_4}) and (\ref{J2_3}), we conclude that
\begin{eqnarray}\label{J1J2}
&& \dbE \int_{Q_b} \cJ_1+\cJ_2dx \nonumber\\
&& \ge  -\dbE \int_{G(t_0-b)}(a+b)^2|\n v(t_0-b,x)|^2dx\nonumber\\
&& \q -\dbE\int_{G(t_0)}a\l e^{C_0}v(t_0,x)^2dx -\dbE \int_{G(t_0-b)}\frac{|x-x_0|^2}{16}v(t_0-b,x)^2 dx\\
&& \q +\dbE \int_{t_0-b}^{t_0} \int_{\pa G(t)} (x-x_0)\cd \nu \[|\n v|^2+\frac{1}{4}(t_0-t+a)e^{2\phi}g_2\]dxdt\nonumber\\
&& \q +\dbE \int_{Q_b} 2\l e^{-C_0} \frac{1}{|\ln (t_0-t+a)|^3} v^2dxdt\nonumber\\
&&\q  -\dbE \int_{Q_b} (t_0-t+a)^2e^{2\phi}\[|\n g_2|^2+\frac{n}{4(t_0-t+a)}g_2^2dt\]dxdt. \nonumber
\end{eqnarray}
\par For $ \cJ_3 $, we have
\begin{eqnarray}\label{J3}
\cJ_3\3n\3n&& = (t_0-t+a)vdv-\(\frac{1}{2}+\frac{n}{4}\)v^2dt-\frac{1}{2}v(x-x_0)\cd\n vdt \nonumber \\
&& = d \[\frac{1}{2}(t_0-t+a)v^2\]+\frac{1}{2}v^2dt-\frac{1}{2}(t_0-t+a)(dv)^2-\(\frac{1}{2}+\frac{n}{4}\)v^2dt  \nonumber\\
&& \q -\n\cd\[\frac{1}{4} (x-x_0)v^2\]dt+\frac{n}{4}v^2dt\\
&& =d \[\frac{1}{2}(t_0-t+a)v^2\]-\n\cd\[\frac{1}{4} (x-x_0)v^2\]dt-\frac{1}{2}(t_0-t+a)(dv)^2. \nonumber
\end{eqnarray}
Consequently, 
\begin{equation}\label{J3_1}
\ba{ll}
\ds \dbE \int_{Q_b} \cJ_3dx \3n&\ds= \frac{a}{2}\dbE\int_{G(t_0)} v(t_0,x)^2dx-\frac{a+b}{2}\dbE \int_{G(t_0-b)}v(t_0-b,x)^2dx \\
\ns&\ds\q -\dbE \int_{Q_b} \frac{1}{2}(t_0-t+a)e^{2\phi} g_2^2dxdt.	
\ea 
\end{equation}
Combining  (\ref{J1J2}) and (\ref{J3}), noting that $ \pa_t\si_a(t)=-\si'(s)|_{s=t_0-t+a} \le 0$, we obtain
\begin{equation}\label{I1I2}
\ba{ll}
\ds \dbE \int_{Q_b} 2\cI_1\cI_2dx=\dbE \int_{Q_b} \cJ_1+\cJ_2+\cJ_3 dx\\
\ns\ds \ge -\cB-\cM+\dbE \int_{Q_b}  \frac{2\l e^{-C_0}}{|\ln (t_0-t+a)|^3}  v^2dxdt.
\ea 
\end{equation}

On the other hand, noting that $ \th(\cd )<2 $,   we find that
\begin{equation*}\label{nu}
\ba{ll}
\ds \frac{1}{2}e^{-2C_0} \dbE \int_{Q_b} \frac{\si_a}{|\ln (t_0-t+a)|^3} |\n u|^2dxdt\\
\ns\ds =\frac{1}{2}e^{-2C_0}\dbE\int_{Q_b} \frac{t_0-t+a}{|\ln (t_0-t+a)|^3}|\n v-\n \phi v|^2dxdt\\
\ns\ds \le e^{-2C_0}\dbE \int_{Q_b} \frac{t_0-t+a}{|\ln (t_0-t+a)|^3}|\n \phi|^2v^2dxdt\\
\ns\ds \q -e^{-2C_0}\dbE \int_{Q_b} \frac{1}{|\ln (t_0-t+a)|^3}v\[-\cI_1-(t_0-t+a)(|\n \phi|^2+\phi_t)v+\frac{1}{2}v\]dxdt\\
\ns\ds \le \dbE \int_{Q_b}\[e^{-2C_0} \frac{1}{|\ln (t_0-t+a)|^3}v\cI_1 +\l e^{-C_0}\frac{1}{|\ln (t_0-t+a)|^3}v^2\]dxdt \\
\ns\ds \le\dbE \int_{Q_b} \[\cI_1^2+\(\l +\frac{1}{4}\)e^{-C_0}\frac{1}{|\ln (t_0-t+a)|^3}v^2\]dxdt.
\ea 
\end{equation*}
This, together with (\ref{p4}) and (\ref{I1I2}),  yields (\ref{Carleman}).
\endpf

\section{Quantitative unique continuation property  for stochastic parabolic equations}

\par In this section, we establish the following quantitative unique continuation results for the equation (\ref{CD}).

\begin{theorem}\label{main_theorem*}
Suppose Assumption  \ref{assumption1} hold. Let $ u \in   L_\dbF^2 (\O; C([0,T]; L^2(G_i(t)))\cap L_\dbF^2(0,T; H^1(G_i(t)))$ be a weak solution of (\ref{CD}). Then there exist  constants $ \eta_1\in (0,1) $ and $ C>0 $, such that for any $ (t_0,x_0)\in G((0,T]) $ and $ r, \rho, R\le \min\{\sqrt{t_0} , R_0, \frac{1}{\sqrt{2e}} \}$ satisfying $ 0<r\le  \rho\le \eta_1 R $ and $ \cB_{R}(x_0) \cap \G= \emptyset $, we have 
\begin{equation}\label{20-4}
	\(\dbE \int_{\cB_\rho(x_0)\cap G(t_0)}u(t_0,x)^2dx\)^{\frac{1}{2}}\le \frac{CR}{\rho}|\ln \rho|^{\frac{3}{2}}\cE^{1-\th}_1\e_1^\th,
\end{equation}
where  $\th=\frac{1}{C(\ln R-\ln r)}$,
\begin{equation}\label{20-3}
	\cE_1\= \(R^{-2}\dbE \int_{t_0-R^2}^{t_0} \int_{\cB_R(x_0)\cap G(t)} u^2dxdt\)^{\frac{1}{2}},\q \e_1\= \(\dbE \int_{\cB_{r}(x_0)\cap G(t_0)} u(t_0,x)^2dx\)^{\frac{1}{2}},
\end{equation}
and  $ \th $ is given in Theorem \ref{main_theorem}.
\end{theorem}

In the proof of Theorem \ref{main_theorem*}, it is necessary to make Assumption \ref{assumption1} when applying the Carleman estimate near the boundary $ I((0,T)) $. However, we can eliminate this assumption if $ (t_0-R^2, t_0)\t \cB_{R}(x_0)\subset G((t_0-R^2,t_0)) $. This leads to a direct and straightforward result as a consequence.

\begin{theorem}\label{main_theorem}
There exist  constants $ \eta_1\in (0,1) $ and $ C> 0 $  such that  for any  $ (t_0,x_0)\in G((0,T)) $ and $ r $, $ \rho $, $ R\le \frac{1}{\sqrt{2e}} $ satisfying $ 0<r\le \rho\le \eta_1 R\le R\le \sqrt t_0 $ and $(t_0-R^2,t_0)\t \cB_R(x_0)\subset G((t_0-R^2,t_0))$, the following estimate holds true:
\begin{equation}\label{2_1}
	\(\dbE\int_{\cB_{\rho}(x_0)} u(t_0,x)^2dx\)^{1/2}\le \frac{CR}{\rho}|\ln \rho|^{\frac{3}{2}}\cE^{1-\th}\e^{\th},
\end{equation}
where
\begin{equation}\label{th1}
	\cE\= \(R^{-2}\dbE \int_{t_0-R^2}^{t_0}\int_{ \cB_R(x_0)}u^2dxdt\)^{1/2},\q \e\=\(\dbE \int_{\cB_r(x_0)} u(t_0, x)^2dx\)^{1/2}.
\end{equation}
\end{theorem}

\begin{remark}
The assumption $ R\le \frac{1}{\sqrt{ 2e}} $ is technical and may be removed. However, it is mentioned here for the sake of convenience in our proof and is sufficient for our applications. On the other hand, the assumption that $R\le \sqrt {t_0}$	
is crucial and cannot be easily removed using the method described in this paper.
\end{remark}

As an immediate corollary of Theorem \ref{main_theorem}, we have the following strong unique continuation property (SUCP for short) for stochastic parabolic equations.
\begin{corollary}\label{SSUCP}
Let $ t_0\in (0,T] $ and $ x_0\in G(t_0) $ such that   $\cB_{r_0}(x_0)\subset G(t_0) $ for some $r_0>0$. 
Let  $ u$ be a weak solution  of (\ref{CD}). If for any $ k\in \dbN $,  there exist a constant $ C_k>0 $ such that
\begin{equation}\label{SSUCP_condition1}
	\dbE \int_{B_{r}(x_0)} u(t_0,x)^2dx\le C_k {r}^{2k}, \q \forall r \in [0,r_0],
\end{equation}
then, $ u(t_0, \cd)\equiv 0 $ in $ G $,  $ \dbP$- a.s. 
\end{corollary}

\begin{remark}
The study of SUCP is  an important subject on partial differential equations and has a long history. However, the stochastic counterpart are far from be well understood. As far as we know, there is no published works addressing this topic.  
\end{remark}

\par {\bf Proof of Corollary \ref{SSUCP} from Theorem \ref{main_theorem}. }Recalling that  for any $N\in\dbN$ and $r>0$, it holds that
\begin{equation}\label{pa-sucp-eq12}
\(\dbE\int_{B_{r}(x_0)}u(t_0,x)^2dxdt\)^{\frac{1}{2}}=O(r^{N}).
\end{equation}
Applying Theorem \ref{main_theorem} to the equation
\eqref{CD},  by \eqref{pa-sucp-eq12},  we
obtain that
\begin{equation}\label{pa-sucp-eq13}
\ba{ll}
\ds \(\dbE\int_{\cB_{\rho}(x_0)}u(t_0,x)^2dx\)^{\frac{1}{2}}\\
\ns\ds \leq \frac{CR}{\rho}\cE^{1-\frac{1}{C(\ln R-\ln r)}}(C_Nr^N)^\frac{1}{C(\ln R-\ln r)}\\
\ns\ds \le \frac{CR}{\rho}\cE^{1-\frac{1}{C(\ln R-\ln r)}}e^{\frac{\ln C_N+N\ln r}{C(\ln R-\ln r)}}, \q
\mbox{for every }N\in\dbN,
\ea 
\end{equation}
Passing (\ref{2_1}) to the limit as $r$ tends to $0$, we get
\begin{equation}\label{20-1}
\(\dbE\int_{\cB_{\rho}(x_0)}u(t_0,x)^2dx\)^{\frac{1}{2}}\le \frac{CR}{\rho}\cE e^{-\frac{1}{C}N}.
\end{equation}
Passing to the
limit as $N\to+\infty$, \eqref{20-1}
yields that $u(t_0,\cd)\equiv 0$ in $\cB_{\rho}$,
$\dbP$-a.s.  By iteration, it follows that $u(t_0, x) \equiv 
0$ in $G$, $\dbP$-a.s.
\endpf

\par Notice that we assume $ t_0\ge R^2 $ in Theorems \ref{main_theorem*} and  \ref{main_theorem}. We also need to deal with the case $  t_0\le R^2$ in our applications. 

\begin{theorem}\label{mth}
Let $ \widetilde G((0,T)) \subset {\bf G}$. Let $ u\in L_\dbF^2 (\O; C([0,T];  L^2(\widetilde G(t)))\cap L_\dbF^2(0,T; H^1(\widetilde G(t))) $solves
\begin{equation}\label{v1}
	\left\{
	\ba{ll}
	\ds du-\D udt=(a_1\cd \n u+b_1 u)dt+c_1u dW(t)  &\mbox{ in } \widetilde G((0,T)),\\
	\ns\ds u(0,x )=0  &\mbox{ on }\widetilde G(0).
	\ea 
	\right.
\end{equation}
Then  there exist  constants $ \eta_1\in (0,1) $ and $ C>0 $, such that for any $ (t_0,x_0)\in \widetilde G((0,T]) $ and $ r, \rho, R\le \min\{ R_0, \frac{1}{\sqrt{2e}} \}$ satisfying $ 0<r\le  \rho\le \eta_1 R $ and $ (\max\{0, t_0-R^2\}, t_0)\t \cB_{R}(x_0) \subset \widetilde G( (\max\{0, t_0-R^2\}, t_0)) $, we have 
\begin{equation}\label{20-6}
	\(\dbE \int_{\cB_\rho(x_0)}u(t_0,x)^2dx\)^{\frac{1}{2}}\le \frac{CR}{\rho}|\ln \rho|^\frac{3}{2}\cE^{1-\th}_2\e_1^\th,
\end{equation}
where
\begin{equation}\label{20-5}
	\cE_2\=\(R^{-2} \dbE \int_{\max\{0, t_0-R^2\}}^{t_0}\int_{\cB_R(x_0)} u^2dxdt\)^{\frac{1}{2}}.
\end{equation}
\end{theorem} 
\begin{remark}
For convenience and  without loss of generality, we make $ \eta_1 $  the same in  Theorems \ref{main_theorem*},   \ref{main_theorem} and   \ref{mth}. 
\end{remark}
\begin{remark}
Let $ u_1$ and $u_2 $ solves (\ref{CD}) in different domains $ G_1((0,T)) $ and $ G_2((0,T)) $, respectively. The reason we need both Theorem \ref{main_theorem*} and Theorem \ref{mth} is that, in the latter sections, we intend to apply Theorem \ref{main_theorem*} to $ u_1$ and $ u_2 $, both away and near the time-varying boundary $ I((0,T)) $, and we intend to apply Theorem \ref{mth} to $ u_1-u_2 $, which solves (\ref{v1}).
\end{remark}
\par Now, we are in the position to  prove of Theorems \ref{main_theorem*} and  \ref{mth}. We introduce the following Lemma first before give the detailed proofs.

\begin{lemma}\label{inital_lemma}

\par (i) 	Under conditions of Theorem \ref{main_theorem*}, for any $ (t_0,x_0)\in G((0,T])$, let $ \rho_2>\rho_1>0 $ satisfying $ \cB_{\rho_2}(x_0)\cap \G=\emptyset $. Then for any $ t\in [\max\{t_0-\frac{R^2}{2}, t_0-\frac{(\rho_2-\rho_1)^2}{18n}\}, t_0] $, we have
\begin{equation}\label{inital_estimate}
	\ba{ll}
	\ds \dbE \int_{\cB_{\rho_1}(x_0)\cap G(t_0)} u(t_0, x)^2dx\\
	\ns\ds \le C\dbE \[\int_{\cB_{\rho_2}(x_0)\cap G(t)}u(t, x)^2dx+\cC_1 \int_{ t_0-R^2}^{t_0}ds\int_{\cB_{\rho_2}( x_0)\cap G(s)} u(s,x)^2dx\],
	\ea 
\end{equation}
where 
$$
\cC_1=\(\frac{1}{R^2}+\frac{1}{(\rho_2-\rho_1)^2}\)\frac{\rho_1^n}{(t_0-t)^{\frac{n}{2}}}e^{-\frac{(\rho_2-\rho_1)^2}{16(t_0-t)}}.
$$
\par (ii) Under conditions of Theorem \ref{mth}, for any $ (t_0, x_0)\in \widetilde G((0,T]) $, let $ \rho_2>\rho_1>0 $ satisfying $ (0,t_0)\t \cB_{\rho_2}(x_0)\subset \widetilde G((0,t_0)) $.  If $ \sqrt{t_0}<\min\{R, \frac{\rho_2-\rho_1}{\sqrt{18n}}\} $, then we have
\begin{equation}\label{inital_estimate1}
	\ba{ll}
	\ds \dbE \int_{\cB_{\rho_1}(x_0)} u(t_0, x)^2dx \le C\cC_2 \int_{ 0 }^{t_0}ds\int_{\cB_{\rho_2}( x_0)} u(s,x)^2dx,
	\ea 
\end{equation}
with 
$$
\cC_2=\frac{\rho_1^n}{t_0^{\frac{n}{2}}(\rho_2-\rho_1)^2}e^{-\frac{(\rho_2-\rho_1)^2}{16t_0}}.
$$
\par (iii)  Under conditions of Theorem \ref{mth}, for any $ (t_0,x_0)\in \widetilde G((0,T]) $, let $ \rho_2>\rho_1>0 $ satisfying $ (0,t_0)\t \cB_{\rho_2}(x_0)\subset \widetilde G((0,t_0)) $. If $  R>\sqrt{t_0}\ge  \frac{\rho_2-\rho_1}{\sqrt{18n}} $, then we have for any $ t\in [t_0-\frac{(\rho_2-\rho_1)^2}{18n},t_0] $,
\begin{equation}\label{20-7}
	\dbE \int_{\cB_{\rho_1}(x_0)} u(t_0,x)^2dx\le C\dbE \[\int_{\cB_{\rho_2}(x_0)}u(t, x)^2dx+\cC_3 \int_{ 0}^{t_0}ds\int_{\cB_{\rho_2}( x_0)} u(s,x)^2dx\],
\end{equation}
with
$$
\cC_3 =\frac{\rho_1^n}{(t_0-t)^{\frac{n}{2}}(\rho_2-\rho_1)^2}e^{-\frac{(\rho_2-\rho_1)^2}{16(t_0-t)}}.
$$
\end{lemma}
\par {\bf Proof of Lemma \ref{inital_lemma}. } Set $ \tilde \rho_1\= \frac{2\rho_1+\rho_2}{3} $ and $ \tilde \rho_2 =\frac{\rho_1+2\rho_2}{3}$. Choose $\f\in  C_0^\infty(\dbR^n)$ such that
\begin{equation}\label{cutoff1}
\f(x)=\left\{\ba{ll}
\ds 1,\q \mbox{in }\cB_{\tilde \rho_1}(x_0),\\
\ns\ds 0, \q \mbox{in } \dbR^n\setminus \cB_{\tilde \rho_2}(x_0).
\ea \right.
\end{equation}
Then $ (\rho_2-\rho_1)|\n \f|+(\rho_2-\rho_1)^2\sum_{j,k=1}^n |\pa_{x_j}\pa_{x_k} \f|  $  is bounded.

Put 
\begin{equation}\label{kernel}
\cK (t,x; t_0, y)\= \frac{1}{(t_0-t)^{n/2}}e^{-\frac{|x-y|^2}{4(t_0-t)}},
\end{equation}
with $ y\in \cB_{\rho_1}(x_0) $. Immediately, we have
\begin{equation}\label{kernel1}
\ba{ll}
\ds \cK_t+\D \cK\3n&\ds =\frac{n}{2}\frac{1}{(t_0-t)^{n/2+1}}e^{-\frac{|x-y|^2}{4(t_0-t)}}-\frac{1}{(t_0-t)^{n/2}}e^{-\frac{|x-y|^2}{4(t_0-t)}}\frac{|x-y|^2}{4(t_0-t)^2}\\
\ns&\ds \q -\frac{1}{(t_0-t)^{n/2}}\n\cd \(\frac{(x-y)}{2(t_0-t)}e^{-\frac{|x-y|^2}{4(t_0-t)}}\) =0.
\ea 
\end{equation}
Let $ v =\f  u $. Then
\begin{equation}\label{equation_u}
dv-\D vdt=(\f a_1\cd \n u+b_1v+2\n \f\cd \n u+\D \f u)dt+c_1vdW(t).
\end{equation}

\par We first deal with case (i), i.e.,   $ R\le \sqrt{t_0} $. 

Let
\begin{equation}\label{25-4-1}
H(t; t_0, y)\= \int_{\cB_{\rho_2}(x_0)\cap G(t)} v(t,x)^2 \cK (t,x; t_0, y)dx, \q y\in \cB_{\rho_1}(x_0).
\end{equation}
From \eqref{25-4-1}, we see that
\begin{equation}\label{29-1}
\ba{ll}
\ds	 \int_{\cB_{\rho_1}(x_0)} H(t; t_0, y)dy\\
\ns\ds= \int_{\cB_{\rho_1}(x_0)}dy \int_{\cB_{\rho_2}(x_0)\cap G(t)} v(t,x)^2\cK (t, x; t_0, y)dx\\
\ns\ds =\int_{\cB_{\rho_2}(x_0)\cap G(t)} v(t,y)^2dy\int_{\cB_{\rho_1}(x_0)} \frac{1}{(t_0-t)^{\frac{n}{2}}}e^{-\frac{|x-y|^2}{4(t_0-t)}}dx.
\ea 
\end{equation}
For almost every fixed sample point $ \omega $, we have $ u\in C([0,T]; L^2(G(t))) $. Then for any $ \cl\e>0 $,  there exist a $ \cl\d_0, \cl\d_1 >0$, for any $ 0<|t_0-t|<\cl\d_0 $, it holds that
\begin{equation}\label{29-3}
\int_{\{\cB_{\rho_1+\cl\d_1}(x_0)\setminus \cB_{\rho_1-\cl\d_1}(x_0)\}\cap G(t)} v^2(t,y)dy<\cl\e,
\end{equation}
\begin{equation}\label{29-2}
\|\int_{\cB_{\rho_1}(x_0)}\frac{1}{(t_0-t)^{n/2}} e^{-\frac{|x-y|^2}{4(t_0-t)}}dx-2^n\pi^{\frac{n}{2}}\|<\cl \e, \mbox{ uniformly when } y\in \cB_{\rho-\cl\d_1}(x_0),
\end{equation}
and that
\begin{equation}\label{29-4}
\sup_{x\in \cB_{\rho_1}(x_0), y\in \dbR^n\setminus \cB_{\rho_1+\cl\d_1}(x_0)}\frac{1}{(t_0-t)^{n/2}}e^{-\frac{|x-y|^2}{4(t_0-t)}}\le \cl\e.
\end{equation}
Let $ \cB_1=\cB_{\rho-\cl\d_1}(x_0) $,
$ \cB_2(t)=\cB_{\rho_1+\cl\d_1}(x_0)\setminus \cB_{\rho_1-\cl\d_1}(x_0) $, and $ \cB_{3}(t)=\cB_{\rho_2}(x_0)\setminus \cB_{\rho_1+\cl \d_1}(x_0) $. It follows from (\ref{29-1}) that for any $ 0<|t_0-t|<\cl\d_0 $,
\begin{eqnarray}\label{29-5}
&& \int_{\cB_{\rho_1}(x_0)} H(t; t_0, y)dy \nonumber \\
&& = \int_{\{\cB_1\cup\cB_2\cup\cB_3\}\cap G(t)}v(t,y)^2dy\int_{\cB_{\rho_1}(x_0)}\frac{1}{(t_0-t)^{\frac{n}{2}}}e^{-\frac{|x-y|^2}{4(t_0-t)}}dx\\
&& \le 2^n\pi^{\frac{n}{2}}\int_{\cB_{\rho_1 }(x_0)\cap G(t)}v(t_0,y)^2dy+C\cl\e. \nonumber
\end{eqnarray}
Similarly, we can get  $$ \int_{\cB_{\rho_1}(x_0)} H(t; t_0, y)dy\ge 2^n\pi^{\frac{n}{2}}\int_{\cB_{\rho_1 }(x_0)\cap G(t)}v(t_0,y)^2dy-C\cl\e.$$
This together with \eqref{29-5}, implies 
\begin{equation}\label{25-4}
\lim_{t\to t_0} \int_{\cB_{\rho_1}(x_0)} H(t; t_0, y)dy = 2^n\pi^{n/2} \int_{\cB_{\rho_1}(x_0)\cap G(t_0)} v(t_0, y)^2dy,\q \dbP\mbox{-}\as
\end{equation}

Set  
$$q_1\=\Vert a_1\Vert_{L^{\infty}_{\dbF}(0,T;C^{n}(O;\dbR^n))}^2 + 2\Vert b_1\Vert_{L^{\infty}_{\dbF}(0,T;C^{n}(O))} + \Vert c_1\Vert_{L^{\infty}_{\dbF}(0,T;C^{n}(O))}^2+ 1.
$$ 
Using integration by parts, we have
\begin{equation}\label{H}
\ba{ll}
\ds d(e^{-q_1s}H(s))\\
\ns \ds =\int_{G(s)}\big[e^{-q_1s}\big(2\cK  vdv +\cK (dv)^2+v^2 \cK_s ds\big)-q_1 e^{q_1s}\cK v^2ds\big]dx\\
\ns\ds =\int_{G(s)} e^{-q_1s}\big[2\cK v(\D v) +\cK (c_1^2+2b_1-q_1)v^2+\cK_sv^2+2\cK v\f a_1\cd \n u\big]dsdx\\
\ns\ds \q + \int_{G(s)} 2e^{-q_1s}\cK v\big(2\n \f\cd \n u+\D \f u\big)dsdx+\int_{G(s)} 2e^{-q_1s}\cK c_1v^2dW(s)dx\\
\ns\ds = \int_{G(s)} e^{-q_1s}\big[-2\cK |\n v|^2+\cK (c_1^2+2b_1-q_1)v^2+2\cK v\f a_1\cd \n u\big]dsdx\\
\ns \ds \q + \int_{G(s)} 2e^{-q_1s}\cK v\big( 2\n \f\cd \n u+\D \f u\big)dsdx+\int_{G(s)} 2e^{-q_1s}\cK c_1u^2dW(s)dx.
\ea 
\end{equation}

Recalling that  $ t\in [\max\{t_0-\frac{R^2}{2}, t_0-\frac{(\rho_2-\rho_1)^2}{18n}\}, t_0] $ and $ v\equiv 0 $ on $  I(t) $,  integrating (\ref{H}) from $ t $ to $ t_0 $, we have that
\begin{equation}\label{20-8}
\ba{ll}
\ds e^{-q_1t_0}2^n\pi^{n/2} \dbE \int_{\cB_{\rho_1}(x_0)\cap G(t_0)}v(t_0, y)^2 dy-e^{-q_1t} \dbE \int_{\cB_{\rho_1}(x_0)}H(t)dy \\
\ns\ds \le  \dbE\int_{\cB_{\rho_1}(x_0)} dy \int_t^{t_0}ds\int_{\{\cB_{\tilde \rho_2}({x_0})\setminus \cB_{\tilde \rho_1}({x_0})\}\cap G(s)} e^{-q_1s}\cK v \big(2\n \f\cd\n u+\D \f u\big)dx.
\ea  
\end{equation} 
Since
$$
(t_0-s)^{-n/2}e^{-\frac{(\rho_2-\rho_1)^2}{36(t_0-s)}}\mbox{ decrease when } s\in [t,t_0]\subset \[t_0-\frac{(\rho_2-\rho_1)^2}{18n}, t_0\],
$$
we get that
\begin{equation}\label{H*1}
\ba{ll}
\ds e^{-q_1 t_0}2^n\pi^{n/2} \dbE \int_{\cB_{\rho_1}(x_0)\cap G(t_0)}v(t_0, y)^2 dy-e^{-q_1 t} \dbE \int_{\cB_{\rho_1}(x_0)}H(t)dy\\
\ns\ds \le C\frac{\rho_1^n}{(t_0-t)^{n/2}}
e^{-\frac{(\rho_2-\rho_1)^{2}}{36(t_0-t)}}
\dbE \int_{t_0-\frac{R^2}{2}}^{t_0}\int_{\cB_{\tilde \rho_2}(x_0)\cap G(s)}\[|\n u|^2+(\rho_2-\rho_1)^{-2}u^2\]dsdx,
\ea 
\end{equation}
which yields that
\begin{equation}\label{H*2}
\ba{ll}
\ds \dbE \int_{\cB_{\rho_1}(x_0)\cap G(t_0)}v(t_0, y)^2 dy\\
\ns\ds \le C\dbE \int_{\cB_{\rho_1}(x_0)}H(t)dy  +C\frac{\rho_1^n}{(t_0-t)^{n/2}}e^{-\frac{(\rho_2-\rho_1)^2}{36(t_0-t)}}\dbE \int_{t_0-\frac{R^2}{2}}^{t_0}\int_{\cB_{\tilde \rho_2}(x_0)\cap G(s)}\[|\n u|^2+(\rho_2-\rho_1)^{-2}u^2\]dsdx.
\ea 
\end{equation}

Let   $ \eta\in  C^\infty_0(\dbR^{n+1}) $ be such that
\begin{equation}\label{eta}
\eta(t,x)=\left\{\ba{ll}
1,\q \mbox{in } (t_0-\frac{R^2}{2},t_0]\t \cB_{ \tilde \rho_2}(x_0),\\
\ns\ds 0,\q  \mbox{in } \dbR^n \setminus \{(t_0-R^2, t_0]\t \cB_{\rho_2}(x_0 )\}.
\ea \right.
\end{equation}
Clearly, $ t_0 \eta_t+(\rho_2-\rho_1)|\n \eta| \le C $ bounded. 

Multiplying (\ref{CD}) by $ \eta^2 u $, we have
\begin{equation}\label{ca1}
\ba{ll}
\ds \eta^2 u(du-\D udt)\\
\ns\ds = d\(\frac{1}{2}\eta^2 u^2\)- \eta\eta_tu^2 -\frac{1}{2}\eta^2(du)^2-\n\cd (\eta^2u\n u)dt+2\eta u \n \eta \cd \n u dt+\eta^2|\n u|^2dt.
\ea 	
\end{equation}
By integrating each side of (\ref{ca1}) on $ (t_0-R^2,t_0)\t \{\cB_{\rho_2}(x_0)\cap G(s)\} $, we obtain
\begin{equation}\label{ca2}
\dbE \int_{t_0-R^2}^{t_0} \int_{\cB_{\rho_2}(x_0)\cap G(s)} \eta^2|\n u|^2dxdt\le  C\[\frac{1}{R^2}+\frac{1}{(\rho_2-\rho_1)^2}\]\dbE \int_{t_0-R^2}^{t_0}\int_{ \cB_{\rho_2}(x_0)\cap G(s) } u^2dxdt, 
\end{equation}
which gives
\begin{equation}\label{caccioppoli}
\dbE \int_{t_0-\frac{R^2}{2}}^{t_0}\int_{\cB_{\tilde \rho_2}(x_0)\cap G(s)}|\n u|^2ds dx\le C\[\frac{1}{R^2}+\frac{1}{(\rho_2-\rho_1)^2}\]\dbE \int_{t_0-R^2}^{t_0} \int_{\cB_{\rho_2}(x_0)\cap G(s)} u^2dsdx.
\end{equation}

Now, by combining (\ref{H*2}) and (\ref{caccioppoli}), we find that
\begin{equation}\label{H2}
\dbE \int_{\cB_{\rho_1}(x_0)\cap G(t_0)}u(t_0, y)^2 dy\le C\[\dbE \int_{\cB_{\rho_1}(x_0)}  H(t)dy+\cC_1 \dbE  \int_{t_0-R^2}^{t_0} ds\int_{\cB_{\rho_2}(x_0)\cap G(s)} u(t,x)^2dx\].
\end{equation}
From \eqref{25-4-1}, we know that
\begin{equation}\label{H3}
\ba{ll}
\ds \dbE \int_{\cB_{\rho_1}(x_0)}  H(t; t_0, y)dy \\ \ns\ds =\dbE \int_{\cB_{\rho_1}(x_0)} dy \int_{\cB_{\rho_2}(x_0)\cap G(t)} v(t,x)^2 \cK (t,x; y)dx\\
\ns\ds =\dbE\int_{\cB_{\rho_2}(x_0)\cap G(t)} v(t,x)^2 dx \int_{\cB_{\rho_1}(x_0)}\frac{1}{(t_0-t)^{n/2}}e^{-\frac{|x-y|^2}{4(t_0-t)}}dy\\
\ns\ds \le C \dbE \int_{\cB_{\rho_2}(x_0)\cap G(t)} u(t,x)^2dx,
\ea 
\end{equation}
Combining (\ref{H2}) and (\ref{H3}), we get (\ref{inital_estimate}).

\ms

\par Next, we handle the  case (ii) that $ \sqrt{t_0}<\min\{R, \frac{\rho_2-\rho_1}{\sqrt{18n}}\} $.  

\ss

Let
$$
\widetilde H(t; t_0, y)\= \int_{\cB_{\rho_2}(x_0) } v(t,x)^2 \cK (t,x; t_0, y)dx, \q y\in \cB_{\rho_1}(x_0).
$$
Similar to (\ref{25-4}), we can obtain
\begin{equation}\label{10-1}
\lim_{t\to t_0} \int_{\cB_{\rho_1}(x_0)} \widetilde H(t; t_0, y)dy = 2^n\pi^{n/2} \int_{\cB_{\rho_1}(x_0)} v(t_0, y)^2dy,\q  a.s.
\end{equation}
Similar to (\ref{H}), we can get
\begin{eqnarray}\label{10-2}
\ds d(e^{-q_1s}\widetilde H (s)) 
\3n&  =& \3n\int_{\widetilde G(s)} e^{-q_1s}\[-2\cK |\n v|^2+\cK (c_1^2+2b_1-c_2)v^2+2\cK v\f a_1\cd \n u\]dsdx \nonumber \\
&&    + \int_{\widetilde G(s)}\! 2e^{-q_1s}\cK v\( 2\n \f\cd \n u\!+\D \f u\)dsdx +\!\int_{G}\! 2e^{-q_1s}\cK c_1u^2dW(s)dx.
\end{eqnarray}
Integrate (\ref{10-2}) from $ 0 $ to $ t_0 $, using the fact  that in this case  function $ (t_0-s)^{-n/2}e^{-\frac{(\rho_2-\rho_1)^2}{36(t_0-s)}} $ decrease on $ [0,t_0] $, noting that $ u(0,\cd)=0 $, we have 
\begin{equation}\label{20-9}
\dbE \int_{\cB_{\rho_1}(x_0)}v(t_0, y)^2 dy\le \frac{C\rho_1^n}{t_0^{n/2}}e^{-\frac{(\rho_2-\rho_1)^2}{36t_0}}\dbE \int_{0}^{t_0}\int_{\cB_{\tilde \rho_2}(x_0)} \[|\n u|^2+(\rho_2-\rho_1)^{-2}u^2\]dsdx.
\end{equation}
Using similar argument as (\ref{caccioppoli}) and noting that $ u(0,\cd)=0 $, we can let the cut-off function $ \eta $ only depend on $ x $, then we can prove the following estimate:
\begin{equation}\label{caccioppoli1}
\dbE \int_0^{t_0}\int_{\cB_{\tilde \rho_2}(x_0)} |\n u|^2dsdx\le \frac{C}{(\rho_2-\rho_1)^{2}}\dbE \int_0^{t_0}\int_{\cB_{\tilde \rho_2}(x_0)}u^2dsdx.
\end{equation}
Combining (\ref{20-9}), (\ref{caccioppoli1}),  we get (\ref{inital_estimate1}).

\ms

Next, we consider the case (iii) that $  R>\sqrt{t_0}\ge  \frac{\rho_2-\rho_1}{\sqrt{18n}} $. 

\ss

Integrating (\ref{10-2}) from $ t $ to $ t_0 $ with $ t\in [t_0-\frac{(\rho_2-\rho_1)^2}{18n}, t_0] $, noting that $ v(0,\cd)\equiv 0 $, we have
\begin{equation}\label{20-11}
\ba{ll}
\ds e^{-q_1t_0}2^n\pi^{n/2} \dbE \int_{\cB_{\rho_1}(x_0)}v(t_0, y)^2 dy-e^{q_1t} \dbE \int_{\cB_{\rho_1}(x_0)}\widetilde  H(t; t_0, y)dy \\
\ns\ds \le \frac{C\rho_1^n}{(t_0-t)^{n/2}}
e^{-\frac{(\rho_2-\rho_1)^{2}}{36(t_0-t)}}
\dbE \int_{0}^{t_0}\int_{\cB_{x_0}(\tilde \rho_2)}\[|\n u|^2+(\rho_2-\rho_1)^{-2 }u^2\]dsdx.
\ea 
\end{equation}
Combining with (\ref{caccioppoli1}) and (\ref{20-11}), we obtain that
\begin{equation}\label{20-12}
\dbE \int_{\cB_{\rho_1}(x_0)}v(t_0, y)^2 dy\le C\[\dbE \int_{\cB_{\rho_1}(x_0)}\widetilde H(t; t_0, y) dy+\cC_3\dbE\int_0^{t_0} ds\int_{\cB_{\rho_2}(x_0)} u(t,x)^2dx\].
\end{equation}
Using similar argument in  (\ref{H3}) to $ \widetilde H $, and noting (\ref{20-12}), we obtain (\ref{20-7}).\endpf 

\ms

\par For a constant $ \tilde \rho>0 $, we set
\begin{equation}\label{D}
D_{\tilde \rho,a}\=\Big\{(t,x)\in (-\i,t_0  )\t \dbR^n\|~ \frac{|x-x_0|^2}{8\l (t_0-t+a)}+\ln (t_0-t+a)< \ln \frac{\tilde \rho^2}{8\l}\Big\}.
\end{equation}
Immediately, for any $(t,x)\in D_{\tilde \rho, a}$, we have
\begin{equation}\label{d1}
\begin{cases}\ds
	\ds |x-x_0|^2< 8\l (t_0-t+a) \ln \frac{\tilde \rho^2}{8\l (t_0-t+a)}\le e^{-1}\tilde \rho^2,\\
	\ns\ds t> t_0+a-\frac{\tilde \rho^2}{8\l}.
\end{cases}
\end{equation}

We first introduce the following lemma.
\begin{lemma}\label{lemma_l}
For $ \phi $ given in (\ref{a1}), we have
\begin{equation}\label{a2}
	e^\phi\ge \(\frac{8\l}{\tilde \rho^2}\)^{\l} , \q \forall (t,x)\in D_{\tilde \rho, a},
\end{equation}
and
\begin{equation}\label{a3}
	(t_0-t+a)^{-k} e^\phi  \le  C^{\l} \(\frac{\l}{\tilde \rho^2}\)^{\l+k},\q \forall (t,x)\in D_{2\tilde\rho,a}\setminus D_{\tilde \rho^2, a},\q k=0, 1, 2.
\end{equation}
\end{lemma}

{\bf Proof of Lemma \ref{lemma_l}. } First, (\ref{a2}) is obvious by noting that $ \si(s)\le s $ given in (\ref{si2}). Now, we focus on the proof of (\ref{a3}).

For $ C_0 $ given in (\ref{si2}), we have
\begin{equation}\label{a4}
e^{\phi}\le e^{C_0 \l} (t_0-t+a)^{-\l}e^{\frac{-|x-x_0|^2}{8(t_0-t+a)}}.
\end{equation}

Let  $\psi(t) = (t_0-t+a)^{-\l-k} e^{-\frac{|x-x_0|^2}{8(t_0-t+a)}} $, $t\in [-R^2,R^2]$. Then $\psi$ attains the maximum value at $ t=t_0+a-\frac{|x-x_0|^2}{8(\l+k)} $. On the other hand,  for any $(t,x)\in  D_{\tilde \rho, a}$ and $|x-x_0|^2=8(\l+k )(t_0-t+a) $, we have
$$
t_0-t+a\le \frac{\tilde \rho^2}{8\l}e^{-\frac{\l+k}{\l}}.
$$
Hence,
\begin{equation}\label{l2}
\ba{ll}
\ds \sup_{ D_{2\tilde \rho,a}\setminus D_{\tilde \rho, a}} (t_0-t+a)^{-k}e^{\phi}\\
\ns\ds \le e^{C_0\l}\sup\Big\{s^{-\l-k}e^{-\l-k}\|~ s\in \Big(\frac{\tilde\rho^2}{8\l e^{\frac{\l+k}{\l}}}, \frac{\tilde\rho^2}{2\l e^{\frac{\l+k}{\l}}} \Big)\Big\}\\
\ns\ds \le e^{C_0\l}\(\frac{8\l}{\tilde \rho^2}\)^{\l+k}e^{\frac{(\l+k)^2}{\l}-\l-k}\\
\ns\ds \le C^{\l} \(\frac{\l}{\tilde \rho^2}\)^{\l+k}.
\ea 
\end{equation}
\endpf

\par {\bf Proof of  Theorem \ref{main_theorem*}.} We divide the proof into four steps.

\ss

{\bf Step 1}. In this step, we give a suitable cut-off function. 

Let  $ R_1\=\frac{1}{2}\sqrt e R $ and  $ r_1<\frac{1}{2}R $ solves $ r_1^2\ln \frac{R_1^2}{r_1^2}=r^2 $. Set $ a=\frac{r_1^2}{8\l} $ and $ b= \frac{R_1^2}{2\l} $,  where $ a$ and $b $ are given in Theorem \ref{theorem_Carleman}. We have $a+b\le \frac{1}{e}$ by choosing $ \l \ge \frac{e}{8}(r_1^2+4R_1^2) $ and $ a<b\le t_0 $ by choosing $ \l\ge \frac{1}{2}+\frac{ e}{4}\ge \frac{1}{2}+\frac{R_1^2}{2t_0} $.

By (\ref{d1}), we obtain
$$
D_{2R_1, a}\subset (t_0-b, t_0)\t \cB_R(x_0)
$$
and
$$
\cB_r(x_0)=\Big\{x\in \dbR^n\|~ (t_0,x)\in D_{R_1, a}\Big\}.
$$

Let 
$$
\psi_{1}(\tau;d_1,d_2) = \begin{cases}
\exp\frac{(d_2-d_1)^2}{(d_1-\tau)(\tau-d_2)}, &\mbox{ if }\tau\in (d_1, d_2),\\
\ns\ds 0 &\mbox{ if }\tau\in \dbR\setminus(d_1, d_2) 	
\end{cases}
$$
and
\begin{equation}\label{psi}
\psi_2(\tau)=\frac{\int_\tau^{d_2}\psi_1(\xi; d_1,d_2)d\xi}{\int_{d_1}^{d_2} \psi_1(\xi,d_1,d_2)d\xi}, \qq  \tau\in \dbR.
\end{equation}
It is easy to check that $ \psi_2(\tau)=1 $ for every $ \tau\in (-\i, d_1) $, $ \psi_2(\tau)=0 $ for every $ \tau\in (d_2, +\i) $. 

Let
\begin{equation}\label{d1d2}
d_1=\ln \frac{(R_1/2)^2}{8\l }, \q d_2=\ln \frac{R_1^2}{8\l }.
\end{equation}
Then, there  exist a constant $ C_1 >0$ such that
\begin{equation}\label{psi1}
\begin{cases}
	\ds |\psi_2'(\tau; d_1,d_2)|\le \frac{C}{d_2-d_1}\le C_1, \q \\
	\ns\ds |\psi_2''(\tau; d_1,d_2)|\le \frac{C}{(d_2-d_1)^2}\le C_1, 
\end{cases}\q \mbox{ for }\tau\in [d_1, d_2].
\end{equation}

Let
\begin{equation}\label{f}
\ba{ll}
\ds \f(t,x)\=\psi_2\(\frac{|x-x_0|^2}{8\l (t_0-t+a)}+\ln (t_0-t+a); d_1,d_2\),\q (t,x)\in (t_0-b, t_0]\t  \cB_R(x_0).
\ea 
\end{equation}
It is easy to check that $ \f\in C_0^\i((t_0-b, t_0]\t  \cB_R(x_0)) $, $ \f=1 $ in $ D_{R_1/2,a } $ and $ \f=0 $ in $ \{(t_0-b, t_0]\t \cB_R(x_0)\}\setminus D_{R_1,a} $. 

Let $ w=\f u $. By  (\ref{CD}), we get
\begin{equation}\label{cutf}
dw-\D wdt= \[\f(a_1\cd\n u+b_1u)+\f_t u-\D \f u-2\n\f\cd\n u\]dt+c_1wdW(t).
\end{equation}

\ss

{\bf Step 2}.  In this step, we apply the Carleman estimate (\ref{Carleman}) to the solution to \eqref{cutf}.

\ss

Applying the Carleman estimate (\ref{Carleman}) to $ w $, and recalling  (ii) of Assumption \ref{assumption1},  we obtain
\begin{equation}\label{c1}
\ba{ll}
\ds \dbE\int_{t_0-b}^{t_0}\int_{\cB_R(x_0)\cap G(t)} \frac{\l e^{-C_0}}{|\ln (t_0-t+a)|^3} e^{2\phi}w^2dxdt\\
\ns\ds \q +\dbE\int_{t_0-b}^{t_0}\int_{\cB_R(x_0)\cap G(t)} \frac{\frac{1}{2}e^{-2C_0}\si_a}{|\ln (t_0-t+a)|^3}e^{2\phi} |\n w|^2dxdt\\
\ns\ds \le e^{2C_0}\dbE \int_{t_0-b}^{t_0}\int_{\cB_R(x_0)\cap G(t)}  \si_a^{2} e^{2\phi} g_1^2dxdt+a\l e^{C_0}\si(a)^{-2\l}\dbE\int_{\cB_r(x_0)\cap G(t_0)} w(t_0,x)^2dx\\
\ds\q +\dbE \int_{t_0-b}^{t_0}\int_{\cB_R(x_0)\cap G(t)} \[\(\frac{1}{2}+\frac{n}{4}\) (t_0-t+a)e^{2\phi}g_2^2 +(t_0-t+a)^2e^{2\phi} |\n g_2|^2\]dxdt, 
\ea 
\end{equation}
where
\begin{equation}\label{g1g2}
\ba{ll}
\ds g_1=\f(a_1\cd\n u+b_1u)+\f_t u-\D \f u-2\n\f\cd\n u,\q g_2=c_1w.
\ea 
\end{equation}
Clearly, 
\begin{equation}\label{21-3}
\f^2|\n u|^2=|\f \n u+\n \f u-\n \f u|^2\le 2|\n w|^2+2|\n \f|^2u^2.
\end{equation}
From (\ref{c1}),  (\ref{21-3}), and noting that
\begin{equation}\label{22-1}
\ba{ll}
\ds \si_a= \sqrt{\si_a} \frac{\sqrt{\si_a}|\ln (t_0-t+a)|^3}{|\ln (t_0-t+a)|^3}\le \frac{CR}{ \l^{\frac{1}{2}}}\frac{1}{|\ln (t_0-t+a)|^3},
\ea 
\end{equation}
we obtain
\begin{eqnarray}\label{c*1}
&& \dbE\int_{t_0-b}^{t_0}\int_{\cB_R(x_0)\cap G(t)} \[\frac{\l}{|\ln (t_0-t+a)|^3} e^{2\phi}w^2+ \frac{\si_ae^{2\phi}}{|\ln (t_0-t+a)|^3} |\n w|^2\]dxdt \nonumber\\
&& \le  \frac{CR}{\l^{\frac{1}{2}}}\int_{t_0-b}^{t_0}\int_{\cB_R(x_0)\cap G(t)}\frac{1}{|\ln (t_0-t+a)|^3} e^{2\phi}(\si_a |\n w|^2+w^2)dxdt\nonumber\\
&& \q +a\l e^{C_0}\si(a)^{-2\l}\dbE\int_{\cB_r(x_0)\cap G(t_0)} w(t_0,x)^2dx\\
&& \q +C\dbE \int_{t_0-b}^{t_0}\int_{\cB_R(x_0)\cap G(t)}(t_0-t+a)^2e^{2\phi}\[(|\n\f|^2+|\D \f|^2+\f_t^2)u^2+|\n \f\cd\n u|^2\]dxdt. \nonumber
\end{eqnarray}

Recalling the definition of $ \f $ in (\ref{f}), we have
\begin{equation}\label{f1}
\ba{ll}
\ds |\n \f|\le \frac{CR}{\l(t_0-t+a)},\\
\ns\ds |\D \f|\le \frac{CR^2}{\l^2 (t_0-t+a)^2}+\frac{C}{\l(t_0-t+a)},\\
\ns\ds |\f_t|\le \frac{CR^2}{\l (t_0-t+a)^2}+\frac{C}{t_0-t+a}.
\ea 
\end{equation}
This, together with (\ref{c*1}), implies that
\begin{equation}\label{c*2}
\ba{ll}
\ds \dbE\int_{t_0-b}^{t_0}\int_{\cB_R(x_0)\cap G(t)} \(\frac{\l}{|\ln (t_0-t+a)|^3} e^{2\phi}w^2+ \frac{\si_a}{|\ln (t_0-t+a)|^3} |\n w|^2\)dxdt\\
\ns\ds \le  \frac{CR}{\l^{\frac{1}{2}}}\int_{t_0-b}^{t_0}\int_{\cB_R(x_0)\cap G(t)}\frac{1}{|\ln (t_0-t+a)|^3} e^{2\phi}(\si_a |\n w|^2+w^2)dxdt\\
\ns\ds \q +a\l e^{C_0}\si(a)^{-2\l}\dbE\int_{\cB_r(x_0)\cap G(t_0)} w(t_0,x)^2dx\\
\ns\ds \q +C\dbE \int_{\{D_{R_1,a}\setminus D_{R_1/2,a}\}\cap  G((t_0-b,t_0))}\Big\{[(t_0-t+a)^{-2}R^4+1]e^{2\phi} u^2+R^2 e^{2\phi} |\n u|^2\Big\}dxdt.
\ea 
\end{equation}
From \eqref{c*2} and Lemma \ref{lemma_l}, we know that there exists a constant $ \l_0 $ such that for $ \l \ge \l_0\ge \max\Big\{\frac{e}{8}(r_1^2+4R_1^2), \frac{1}{2}+\frac{ e}{4}\Big\} $,  it holds that
\begin{eqnarray}\label{c2}
&& \dbE\int_{t_0-b}^{t_0}\int_{\cB_R(x_0)\cap G(t)} \frac{\l}{|\ln (t_0-t+a)|^3} e^{2\phi}w^2dxdt \nonumber\\
&& \le C r_1^2\si(a)^{-2\l}\dbE\int_{\cB_r(x_0)\cap G(t_0)} u(t_0,x)^2dx \nonumber\\
&& \q +C^{\l} \(\frac{\l}{R^2}\)^{2\l}\dbE \int_{\{D_{R_1,a}\setminus D_{R_1/2,a}\}\cap G((t_0-b,t_0))}\big(R^2|\n u|^2+u^2\big)dxdt\\
&& \le C^{\l} \frac{\l^{2\l}}{R^{4\l}}\dbE \int_{\{D_{R_1,a}\setminus D_{R_1/2,a}\}\cap G((t_0-b,t_0))}\big(R^2|\n u|^2+u^2\big)dxdt \nonumber\\
&& \q 
+C^{\l}\frac{\l^{2\l}}{r_1^{4\l-2}}\dbE \int_{\cB_r(x_0)\cap G(t_0)} u(t_0,x)^2dx.\nonumber
\end{eqnarray}
and for any $ \tilde \rho\le \frac{R_1}{2} $, 
\begin{eqnarray}\label{c3}
&&\dbE\int_{t_0-b}^{t_0}\int_{\cB_R(x_0)\cap G(t)} \frac{\l}{|\ln (t_0-t+a)|^3} e^{2\phi}w^2dxdt \nonumber \\
&& \ge \dbE \int_{D_{\tilde \rho,a}\cap G ((t_0-b,t_0))} \frac{\l}{|\ln (t_0-t+a)|^3}e^{2\phi}u^2dxdt\\
&& \ge \frac{1}{C^\l}\frac{\l^{2\l}}{\tilde \rho^{4\l}} \dbE \int_{D_{\tilde \rho,a}\cap G((t_0-b,t_0))} \frac{1}{ |\ln (t_0-t+a)|^3} u^2dxdt. \nonumber 
\end{eqnarray}

\ss

{\bf Step 3}. In this step, we prove the following inequality:
\begin{equation}\label{12-1}
\ba{ll}
\ds \dbE \int_{D_{R_1,a}\cap  G((t_0-R^2, t_0))}|\n u|^2dxdt+ \dbE\int_{\cB_{R/2}(x_0)\cap G(t_0)}u(t_0,x)^2dx\\
\ns\ds \le \frac{C}{R^2} \dbE \int_{t_0-R^2}^{t_0}\int_{\cB_R(x_0)\cap G(t)} u^2dxdt,\q \forall \l \ge \l_0.
\ea 
\end{equation}

By (\ref{d1}), we have  
\begin{equation}\label{21-1}
D_{R_1,a}\subset \(t_0-\frac{R^2}{2}, t_0\)\t \cB_{\frac{R}{2}} (x_0).
\end{equation}
Let $ \tilde\f \in C^\i((t_0-R^2, t_0]\t \cB_R(x_0)) $ such that
\begin{equation}\label{21-2}
\tilde \f(t,x)=\left\{\ba{ll}
\ds 1,& (t,x)\in  (t_0-\frac{R^2}{2}, t_0]\t \cB_{\frac{R}{2}} (x_0),\\
\ns\ds 0,&(t,x) \in  (t_0-R^2, t_0]\t \cB_R(x_0)\setminus (t_0-3R^2/4,t_0]\t \cB_{3R/4}(x_0).
\ea\right.
\end{equation}
Then $ R^2|\tilde \f_t|+R|\n \tilde \f| $ is bounded.

Similar to the proof of (\ref{caccioppoli}), by multiplying (\ref{CD}) by $ \tilde \f^2 u $,   we get that
\begin{equation}\label{12-2}
\ba{ll}
\ds \tilde\f^2(du-\D udt)\\
\ns\ds =d\(\frac{1}{2} \tilde\f^2 u^2\)-\tilde \f\tilde \f_t u^2-\frac{1}{2}\tilde \f^2 (du)^2-\n\cd (\tilde\f^2 u\n u)dt+2\tilde \f u\n\tilde \f\cd\n udt+\tilde\f^2|\n u|^2dt.
\ea 
\end{equation}
Integrating each side of (\ref{12-2}) on $ (t_0-R^2, t_0)\t \{\cB_R(x_0)\cap G(t)\} $,  we obtain (\ref{12-1}). 

\ss

{\bf Step 4}.  In this step, we complete the proof.

\ss

Combining (\ref{c2}), (\ref{c3}) and \eqref{12-1},  we obtain that there exist a constant $ \l_0^*>2\l_0 $ such that for any $ \l\ge \l_0^* $, it holds that
\begin{equation}\label{c4}
\ba{ll}
\ds \dbE \int_{D_{\tilde  \rho, a}} \frac{1}{ |\ln (t_0-t+a)|^3} u^2dxdt\\
\ns\ds \le  \(\frac{C\tilde \rho}{R}\)^{2\l}\dbE \int_{\{D_{R_1,a}\setminus D_{R_1/2,a}\}\cap G((t_0-b,t_0))} u^2dxdt
+r_1^2\(\frac{C\tilde \rho}{r_1}\)^{2\l}\dbE \int_{\cB_r(x_0)\cap  G(t_0)} u(t_0,x)^2dx.
\ea 
\end{equation}

Let
\begin{equation}\label{t1}
\ba{ll}
\ds t_1\=t_0-\frac{\tilde \rho^2 e^{-1/2}}{8\l }+a,\q t_2\=t_0-\frac{\tilde \rho^2e^{-1}}{8\l }+a,\\
\ns\ds \cQ \=\{(t_1, t_2)\t\cB_{\frac{\tilde \rho}{2\sqrt[4] e}}(x_0)\}\cap G((t_1, t_2)).
\ea 
\end{equation}
It is easy to check that $ 	\cQ \subset D_{\tilde \rho, a} $.
Moreover, by  setting $ \eta_1<\frac{1}{e} $, we obtain that
$$
\frac{1}{e}r^2\(\ln \frac{R_1^2}{r^2}\)^{-1}\le r_1^2\le r^2\(\ln \frac{R_1^2}{r^2}\)^{-1}\le  r^2.
$$
From (\ref{c4}), recalling the definition of $ \cE_1$ and $\e_1 $ in (\ref{20-3}), we get 
\begin{equation}\label{to1}
\( \dbE \int_\cQ u^2dxdt\)^{1/2}\le  |\ln\tilde \rho|^{\frac{3}{2}}\[R\(\frac{C\tilde \rho}{R}\)^{\l}\cE_1 +r\(\frac{C\tilde\rho}{r(
	\ln R-\ln r)^{-1/2}}\)^{\l}\e_1\].
\end{equation}

By applying Lemma \ref{inital_lemma} to the left hand side of (\ref{to1}) and  let $ \rho=\frac{\tilde \rho}{4\sqrt[4] e} $, we get that there exists a $ \l_1>\l_0^* $, such that for any $ \l\ge \l_1 $ and  $ t\in (t_1, t_2) $, it holds that
\begin{equation}\label{to2}
\dbE \int_{\cB_{\frac{\tilde \rho}{2\sqrt[4] e}}(x_0)\cap  G(t)}u(t,x)^2dx\ge \frac{1}{C}\dbE \int_{\cB_{\rho}(x_0)\cap  G(t_0)} u(t_0,x)^2dx-\(\frac{C}{R^2}+\frac{C}{\rho^2}\)\frac{R^2\rho^n}{(t_0-t)^{\frac{n}{2}}}e^{-\frac{C\rho^2}{(t_0-t)}}\cE_1^2.
\end{equation}
From (\ref{t1}) and recalling $ a=\frac{r_1^2}{8\l } $, we obtain that there exist a constant $ \l_2\ge \l_1 $, such that for any $ \l\ge \l_2 $, 
\begin{equation}\label{to*3}
\dbE \int_{\cB_{\frac{\tilde \rho}{2\sqrt[4] e}}(x_0)\cap G(t)}u(t,x)^2dx\ge \frac{1}{C}\dbE \int_{\cB_{\rho}(x_0)\cap G(t_0)} u(t_0,x)^2dx-\frac{CR^2}{\rho^2}e^{-C\l}\cE_1^2.
\end{equation}
By integrating $ t $ from $ t_1 $ to $ t_2 $,  we find that
\begin{equation}\label{to6}
\dbE \int_{\cQ} u^2dxdt\ge \frac{\rho^2}{\l }\[\frac{1}{C}\dbE \int_{\cB_{\rho}(x_0)\cap G(t_0)} u(t_0,x)^2dx-\frac{CR^2}{\rho^2}e^{-C\l}\cE_1^2\].
\end{equation}
Combining (\ref{to6}) with (\ref{to1}), we conclude that 
\begin{equation}\label{to3}
\(\dbE \int_{\cB_{\rho}(x_0)\cap G(t_0) } u(t_0,x)^2dx\)^{1/2}\le Ce^{-C\l}\frac{R}{\rho}\cE_1 +|\ln \rho|^{\frac{3}{2}}\[ \(\frac{C\rho}{R}\)^{\l}\cE_1 +\(\frac{C\rho}{r(\ln R-\ln r)}\)^{\l}\e_1\].
\end{equation}
Now, we choose $ \eta_1 <\frac{1}{e} $ small enough such that $ C\rho\le C\eta_1R<R  $ and  set
$$
\l_2=\frac{\ln \cE_1 -\ln \e_1}{\ln R-\ln (r(\ln R-\ln r)^{-\frac{1}{2}})}.
$$

If $ \l_2\ge\l_1  $, choosing $ \l=\l_2 $,   by (\ref{to3}), we get
\begin{equation}\label{to4}
\(\dbE\int_{\cB_{\rho}(x_0)\cap G(t_0)} u(t_0,x)^2dx\)^{1/2}\le |\ln \rho|^{\frac{3}{2}} \cE_1^{1-\th_1}\e_1^{\th_1}+\frac{CR}{\rho}\cE_1^{1-\th_2}\e_1^{\th_2},
\end{equation}
where $$ \th_1=\frac{\ln R-\ln C\rho}{\ln R-\ln (r(\ln R-\ln r)^{-\frac{1}{2}})}, \qq \th_2=\frac{C}{\ln R-\ln (r(\ln R-\ln r)^{-\frac{1}{2}})}. $$ By \eqref{12-1}, we have 
\begin{equation}\label{to4-1}
\e_1\le C\cE_1.
\end{equation}
Notice that there exist constant $ C > 0 $ such that
\begin{equation}\label{th3}
\th_1\ge \th=\frac{1}{C(\ln R-\ln r)}, \qq \th_2 \ge \th=\frac{1}{C(\ln R-\ln r)}.
\end{equation}
This, together with \eqref{to4} and \eqref{to4-1}, yields
\begin{equation}\label{to7}
\(\dbE\int_{\cB_{\rho}(x_0)\cap G(t_0)} u(t_0,x)^2dx\)^{1/2}\le \frac{CR}{\rho}|\ln \rho|^{\frac{3}{2}} \cE_1^{1-\th}\e_1^{\th}.
\end{equation}

If $ \l_2< \l_1   $, by \eqref{12-1} again,  we have
\begin{equation}\label{to5}
\ba{ll}
\ds \(\dbE\int_{\cB_{\rho}(x_0)\cap G(t_0)} u(t_0,x)^2dx\)^{1/2}\\
\ns\ds \le C\cE_1^{1-\th_1}\cE_1^{\th_1}\le C\cE_1^{1-\th_1}\(e^{C(\ln R-\ln (r(\ln R-\ln r)))}\e_1\)^{\th_1}\le C\cE_1^{1-\th}\e_1^{\th}.
\ea 
\end{equation}
Combining (\ref{to7}) and (\ref{to5}), we obtain (\ref{20-4}). \endpf

\ms

\par {\bf Proof of  Theorem \ref{mth}. } When $ t_0\ge R^2 $, the proof is exactly the same as Theorem \ref{main_theorem*}. We only need to focus on the case that $ t_0<R^2<\frac{1}{2e} $. We still let $ R_1=\frac{1}{2}\sqrt e R $ and $ r_1<\frac{1}{2} R $ solves $ r_1^2\ln \frac{R_1^2}{r_1^2}=r^2 $. Set $ a=\frac{r_1^2}{8\l} $ and $ b=t_0 $, where $ a$ and $b $ are given in Theorem \ref{theorem_Carleman}. By choosing $ \l >\frac{e}{4} r_1^2 $, we get $ a+b\le \frac{1}{e} $. Similarly, it is easy to check that
\begin{equation}\label{22-2}
D_{2R_1, a}\subset (-\i, t_0)\t \cB_R(x_0),
\end{equation}
and
\begin{equation}\label{22-3}
\cB_r(x_0)=\{x\in \dbR^n|~ (t_0,x)\in D_{R_1,a}\}.
\end{equation}

Let $ \f$ and $ w $ be given by (\ref{f}) and (\ref{cutf}), respectively. Applying the Carleman estimate (\ref{Carleman}) to $ w $, recalling Assumption \ref{assumption1} and   $ u(0,\cd)=0 $, we obtain that
\begin{equation}\label{22-4}
\ba{ll}
\ds \dbE\int_{0}^{t_0}\int_{\cB_R(x_0)} \frac{\l e^{-C_0}}{|\ln (t_0-t+a)|^3} e^{2\phi}w^2dxdt +\dbE\int_{0}^{t_0}\int_{\cB_R(x_0)} \frac{\frac{1}{2}e^{-2C_0}\si_a}{|\ln (t_0-t+a)|^3}e^{2\phi} |\n w|^2dxdt\\
\ns\ds \le e^{2C_0}\dbE \int_{0}^{t_0}\int_{\cB_R(x_0)}  \si_a^{2} e^{2\phi} g_1^2dxdt+a\l e^{C_0}\si(a)^{-2\l}\dbE\int_{\cB_r(x_0)} w(t_0,x)^2dx\\
\ds\q +\dbE \int_{0}^{t_0}\int_{\cB_R(x_0)} \[\(\frac{1}{2}+\frac{n}{4}\) (t_0-t+a)e^{2\phi}g_2^2 +(t_0-t+a)^2e^{2\phi} |\n g_2|^2\]dxdt, 
\ea 
\end{equation}
with $ g_1$ and $g_2 $ are given in (\ref{g1g2}).

For any $ \tilde \rho>0 $, set
\begin{equation}\label{22-6}
\widetilde D_{\tilde \rho, a}=D_{\tilde \rho, a}\cap \widetilde G((0,t_0)).
\end{equation}
Using the same argument from (\ref{c*1})  to (\ref{c3}), we can get that there exists a constant $ \l_0\ge \frac{er_1^2}{4} $ such that for any $ \l\ge \l_0 $, it holds that
\begin{equation}\label{22-5}
\ba{ll}
\ds \dbE \int_0^{t_0}\int_{\cB_R(x_0)} \frac{\l}{|\ln (t_0-t+a)|^3}e^{2\phi}w^2dxdt\\
\ns\ds \le C^\l \frac{\l^{2\l}}{R^{4\l}}\dbE \int_{\widetilde D_{R_1,a}\setminus \widetilde D_{R_1, a}} \big(R^2|\n u|^2+u^2\big)dxdt+ C^\l\frac{\l^{2\l}}{r_1^{4\l-2}}\dbE\int_{\cB_r(x_0)} u(t_0,x)^2dx,
\ea 
\end{equation}
and for any $ \tilde \rho\le\frac{R_1}{2} $,
\begin{eqnarray}\label{22-7}
&&\dbE \int_0^{t_0}\int_{\cB_R(x_0)} \frac{\l}{|\ln (t_0-t+a)|^3}e^{2\phi}w^2dxdt \nonumber\\
&& \ge \dbE \int_{\widetilde D_{\tilde \rho, a}}\frac{\l}{|\ln(t_0-t+a)|^3}e^{2\phi}u^2dxdt\\
&& \ge \frac{1}{C^\l}\frac{\l^{2\l}}{\tilde \rho^{4\l}}\dbE \int_{\widetilde  D_{\tilde \rho, a}} \frac{1}{|\ln (t_0-t+a)|^3}u^2dxdt. \nonumber
\end{eqnarray}

Replacing $ \tilde \f $ in (\ref{21-2}) by a $\hat \f\in C_0^\infty(\dbR^n)$ satisfying
$$
\hat \f(x)=\left\{\ba{ll}
\ds 1, &\mbox{ in }\cB_{\frac{R}{2}}(x_0),\\
\ns\ds 0, &\mbox{ in }\dbR^n\setminus \cB_R(x_0),
\ea \right.
$$
repeating exactly the same procedure in the proof of \eqref{12-1}, we can get the following inequality:
\begin{equation}\label{22-8}
\ba{ll}
\ds \dbE \int_{\widetilde D_{R_1,a}}|\n u|^2dxdt+ \dbE\int_{\cB_{R/2}(x_0)}u(t_0,x)^2dx \le \frac{C}{R^2} \dbE \int_{0}^{t_0}\int_{\cB_R(x_0)} u^2dxdt.
\ea 
\end{equation}
Combining (\ref{22-8}) with (\ref{22-5}) and (\ref{22-7}), we conclude that for any $ \l\ge 2\l_0 $, 
\begin{equation}\label{22-9}
\ba{ll}
\ds \dbE \int_{\widetilde D_{\tilde  \rho, a}} \frac{1}{ |\ln (t_0-t+a)|^3} u^2dxdt\\
\ns\ds \le  \(\frac{C\tilde \rho}{R}\)^{2\l}\dbE \int_{\widetilde D_{R_1,a}\setminus \widetilde D_{R_1/2,a}} u^2dxdt
+r_1^2\(\frac{C\tilde \rho}{r_1}\)^{2\l}\dbE \int_{\cB_r(x_0)} u(t_0,x)^2dx.
\ea 
\end{equation}

Let us use the same notations as (\ref{t1}), i.e.,
\begin{equation}\label{t1*}
\ba{ll}
\ds t_1=t_0-\frac{\tilde \rho^2 e^{-1/2}}{8\l }+a,\q t_2=t_0-\frac{\tilde \rho^2e^{-1}}{8\l }+a,\q
\cQ_1 \=(t_1, t_2)\t\cB_{\frac{\tilde \rho}{2\sqrt[4] e}}(x_0).
\ea 
\end{equation}

In what follows, we deal with  two different cases. 

\ss

If $ t_0\le \frac{\tilde \rho^2 e^{-1/2}}{8\l } $, by choosing $ \tilde \rho=\rho $ and $ \l\ge \frac{9ne^{-\frac{1}{2}}}{4(1-\frac{\sqrt e}{4})^2} $, it is easy to check that $ \sqrt t_0\le \min\{R, \frac{R-\rho}{\sqrt{18n}}\}  $. Then we can apply (ii) in Lemma \ref{inital_lemma} to get
\begin{equation}\label{23-1}
\ba{ll}
\ds \dbE \int_{\cB_{\rho}(x_0)} u(t_0,x)^2dx\\
\ns\ds \le C \frac{\l^{\frac{n}{2}}e^{-C\l}}{R^2}\int_0^{t_0}ds\int_{\cB_R(x_0)} u(s,x)^2dx\\
\ns\ds \le Ce^{-C\l} \cE_2^2.
\ea 
\end{equation}

\ss

If $ t_0> \frac{\tilde \rho^2 e^{-1/2}}{8\l }$, it is easy to check that $ \cQ_1\subset \widetilde D_{\tilde\rho, a} $. From (\ref{22-9}), we have
\begin{equation}\label{23-2}
\( \dbE \int_{\cQ_1} u^2dxdt\)^{1/2}\le  |\ln\tilde \rho|^{\frac{3}{2}}\[R\(\frac{C\tilde \rho}{R}\)^{\l}\cE_2 +r\(\frac{C\tilde\rho}{r(
	\ln R-\ln r)^{-1/2}}\)^{\l}\e_1\]. 
\end{equation}
By choosing $ \rho=\frac{\tilde \rho}{4\sqrt[4] e} $, using similar argument like (\ref{to2}) to (\ref{to6}),  applying (iii) in Lemma \ref{inital_lemma}, we obtain
\begin{equation}\label{23-3}
\dbE \int_{\cQ_1} u^2dxdt\ge \frac{\rho^2}{\l }\(\frac{1}{C}\dbE \int_{\cB_{\rho}(x_0)} u(t_0,x)^2dx-\frac{C}{\rho^2}e^{-C\l}\cE_2^2\)
\end{equation}

Combining (\ref{23-1}), (\ref{23-2}) and (\ref{23-3}), we conclude that
\begin{equation}\label{2304}
\(\dbE \int_{\cB_{\rho}(x_0) } u(t_0,x)^2dx\)^{1/2}\le Ce^{-C\l}\cE_2 +|\ln \rho|^{\frac{3}{2}}\[ \(\frac{C\rho}{R}\)^{\l}\cE_2 +\(\frac{C\rho}{r(\ln R-\ln r)^{-\frac{1}{2}}}\)^{\l}\e_1\].
\end{equation}
Using exactly the same argument like (\ref{to4}) to (\ref{to5}), we obtain (\ref{20-6}). \endpf

\par To the end of this section, we prove the following so-called small propagation estimate which will be useful later.

\begin{definition}
Let $ G \subset \dbR^n$ be a bounded domain. We call $ \pa G $ is of Lipschitz class with constant $ \rho_0, \a $ if for any $ x_0\in \pa G $, there exist $ \zeta\=\zeta(x_0)\in \dbR^n $, $ |\zeta|=1 $ such that
\begin{equation}\label{14-1}
	\Big\{x\in \cB_{\rho_0}(x_0)\|~ \frac{(x-x_0)\cd \zeta}{|x-x_0|}> \cos\a \Big\}\subset G. 
\end{equation}
\end{definition}
\begin{theorem}\label{theoremC1}
Suppose Assumptions \ref{assumption2}--\ref{assumption3} hold.  For $i=1,2$, let  $ u_i $ be the solution of (\ref{CD}) with the unknown boundary $ I_i(t) $. Let $ t_0\in (0,T] $, $ x_0\in \pa G_i(t_0) $ and suppose $ \pa G_i(t_0) $ is of Lipschitz class with constant $ \rho_0, \a$ satisfying $ \sqrt t_0-\frac{\rho_0\sin\a}{E(1+\sin \a)}\ge 0 $. Set
\begin{equation}\label{14-3}
	\mu_1\=\frac{\rho_0}{1+\sin \a},\q w_1\= x_0+\mu_1\zeta, \q \rho_1\= \frac{1}{4E}\eta_1 \mu_1 \sin \a, \q  \si\= \dbE \int_{\cB_{\rho_1}(w_1)} u_i(t_0, x)^2 dx,
\end{equation}
where $\eta_1$ is the constant in Theorem \ref{main_theorem*} and $ E $ is the constant in Assumption \ref{assumption3}. If $ \si\in (0,1) $, then we have
\begin{equation}\label{14-2}
	\dbE |u_i(t_0, x_0)|^2\le C\k_0^C|\ln \si |^{-\frac{\a}{C}}.
\end{equation}
\end{theorem}
\begin{remark}
Here we assume that $ \sqrt{t_0}-\frac{\rho_0\sin\a}{E(1+\sin \a)}\ge 0 $. The reason is that we need $ t_0-4\eta_1^{-1}\rho_1\ge 0 $ when we apply Theorem \ref{main_theorem*}. Clearly, from Assumption \ref{assumption1}, we know that $ \pa G_i(t_0) $ is of Lipschitz class with the constants $ R_0$ and $\frac{\pi}{4} $.  Then, for any  $  \rho_0<R_0$, $  \a<\frac{\pi}{4} $,   $ \pa G_i(t_0) $ is of Lipschitz class with  the constants $ \rho_0$ and  $\a $ immediately from the definition. Thus,  we can choose $ \rho_0 $ and $ \a $ depending on $ t_0 $ and small to make the assumptions hold.
\end{remark}

We introduce the following lemma we will be useful later.
\begin{lemma}\label{lemx}
If a positive sequence $ \{x_n\}_{n=1}^{+\i} $ satisfy the recursive  
relation $ x_n\le C_1 x_{n-1}^s $, $ n>1  $ for some $ C_1>1 $ and 
$ s<1 $, then we have
\begin{equation}\label{b1}
	x_n\le C_1^{\frac{1}{1-s}}x_1^{s^{n-1}}.
\end{equation}
\end{lemma}

{\bf Proof of Lemma \ref{lemx}. } We have
$$
\ln x_n-\frac{\ln C_1}{1-s}\le s\(\ln x_{n-1}-\frac{\ln C_1}{1-s}\).
$$
Then, 
$$
\ln x_n\le \frac{\ln C_1}{1-s}+\(\ln x_1-\frac{\ln C_1}{1-s}\)s^{n-1},
$$
which yields
$$
x_n\le e^{\frac{\ln C_1}{1-s}}e^{-\frac{\ln C_1}{1-s}s^{n-1}} 
x_{1}^{s^{n-1}}\le C_1^{\frac{1}{1-s}}x_1^{s^{n-1}}.
$$
\endpf

{\bf Proof of Theorem \ref{theoremC1}. } The main idea to prove (\ref{14-2}) is  to irritate the  two sphere  one cylinder inequality given in Theorem \ref{main_theorem*} inside the conic domain like (\ref{14-1}). Without loss of generality, we assume that $ x_0=0 $.  It is easy to check that
\begin{equation}\label{14-4}
\cB_{4\eta_1^{-1}E\rho_1}(w_1)\subset \Big\{x\in \cB_{\rho_0}(x_0)\|~ \frac{x\cd \zeta}{|x-x_0|}> \cos\a \Big\}\subset G_i(t),\q i=1,2.
\end{equation}
From Assumption \ref{assumption3} and $ t_0-4\eta_1^{-1}\rho_1\ge 0 $ , we have
\begin{equation}\label{15-2}
\( t_0-4\eta_1^{-1}\rho_1, t_0\)\t \cB_{4\eta_1^{-1}\rho_1} (w_1)\subset G_i\(( t_0-4\eta_1^{-1}\rho_1, t_0)\).
\end{equation}
Applying the inequality (\ref{20-3}) with $ R=4\eta_1^{-1}\rho_1, ~\rho= 3\rho_1,~ r=\rho_1 $, we obtain that
\begin{equation}\label{15-3}
\ba{ll}
\ds \dbE \int_{\cB_{3\rho_1}(w_1)} u_i(t_0,x)^2dx \le C\k_0^{2-s}\( \dbE\int_{\cB_{\rho_1}(w_1)}u_i(t_0,x)^2dx\)^s,
\ea 
\end{equation}
where $ \k_0 $ appears in Assumption \ref{assumption2} and $ s=\frac{C}{\ln 4\eta_1^{-1}} $ is a  constant. 

Set
\begin{equation}\label{14-5}
a\=\frac{1-\frac{1}{4E}\eta_1\sin\a}{1+\frac{1}{4E}\eta_1 \sin\a },\q \mu_2\=a \mu_1, \q w_2\= \mu_2\zeta, \q \rho_2\=a\rho_1.
\end{equation}
It is easy to check that
\begin{equation}\label{14-6}
\cB_{\rho_2}(w_2)\subset \cB_{3\rho_1}(w_1)\subset \cB_{4\eta_1^{-1}\rho_1}(w_1).
\end{equation}
From (\ref{15-3}), we obtain
\begin{equation}\label{14-7}
\dbE \int_{\cB_{\rho_2}(w_2)} u_i(t_0,x)^2dx\le C\k_0^2\( \dbE\int_{\cB_{\rho_1}(w_1)}u_i(t_0,x)^2dx\)^s.
\end{equation}

For $ k\in \dbN $, put
\begin{equation}\label{14-8}
\mu_k=a^{k-1}\mu_1,\q w_k=\mu_k\zeta,\q \rho_k=a^{k-1}\rho_1.
\end{equation}
It is easy to check that
\begin{equation}\label{14-8*}
\cB_{\rho_{k+1}}(w_{k+1})\subset \cB_{3\rho_k}(w_k)\subset \cB_{4\eta_1^{-1}\rho_k}(w_k)\subset\cB_{4\eta_1^{-1}E\rho_k}(w_k) \subset\Big\{x\in \cB_{\rho_0}(x_0)\|~ \frac{(x-x_0)\cd \zeta}{|x-x_0|}> \cos\a \Big\}.
\end{equation}
Similar to \eqref{15-3}, we get
\begin{equation}\label{14-9}
\dbE \int_{\cB_{\rho_{k+1}}(w_{k+1})} u_i(t_0,x)^2dx\le C\k_0^2\( \dbE\int_{\cB_{\rho_k}(w_k)}u_i(t_0,x)^2dx\)^s.
\end{equation}
For a real number $ \tilde \si<\mu_1-\rho_1 $, let $ \overline k $ be the smallest positive integer such that $ \mu_k-\rho_k< \tilde \si $. Then We have
\begin{equation}\label{14-10}
\frac{\ln \tilde \si-\ln (\mu_1-\rho_1)}{\ln a}+1\le\overline k\le   \frac{\ln \tilde \si-\ln (\mu_1-\rho_1)}{\ln a}+2.
\end{equation}
Irritating the inequality (\ref{14-9})  $ \overline k$ times, by Lemma \ref{lemx}, we find that   there exist a constant $ C>0 $ such that
\begin{equation}\label{14-11}
\dbE \int_{\cB_{\rho_{\overline k}}(w_{\overline k})}u_i(t_0,x)^2dx\le C\k_0^{\frac{2}{1-s}}\( \dbE\int_{\cB_{\rho_1}(w_1)}u_i(t_0,x)^2dx\)^{(\frac{\tilde \si}{\mu_1-\rho_1})^\frac{\ln s}{\ln a}}.
\end{equation}
{Recall that $u_i\in C_\dbF([0,T]; L^2(\O; C^1(G_i(t)))$, then we have for any $ x\in \cB_{\rho_{\overline k}}(w_k) $, 
\begin{equation}\label{14-12}
	\ba{ll}
	\ds \dbE \Vert u_i(t_0,\cd)\Vert_{L^\i(\cB_{\rho_{\overline k}}(w_k))}^2\\
	\ns\ds \le \dbE \Vert u_i(t_0,\cd)-u_i(t_0, x)\Vert_{L^\i(\cB_{\rho_{\overline k}}(w_k))}^2+\dbE u_i(t_0,x)^2\\
	\ns\ds \le  \dbE \Vert \n u_i\Vert_{L^\i(\cB_{\rho_{\overline k}}(w_k))}^2 \rho_{\cl k}^2+\dbE u_i(t_0,x)^2\\
	\ns\ds \le  C{\k_0^2}\rho_{\cl k}^2+\dbE u_i(t_0,x)^2.
	\ea 
\end{equation}
Combing (\ref{14-11}) with  (\ref{14-12}), we obtain
\begin{eqnarray}\label{14-13}
	&&\dbE \Vert u_i(t,\cd)\Vert_{L^\i(\cB_{\rho_{\overline k}}(w_k))}^2 \nonumber\\
	&& \le C\k_0^{\frac{2}{1-s}}\( \dbE\int_{\cB_{\rho_1}(w_1)}u_i(t_0,x)^2dx\)^{(\frac{\tilde \si}{\mu_1-\rho_1})^\frac{\ln s}{\ln a}}+C\k_0^2\rho_{\overline k}^2\\
	&&\le C\k_0^{\frac{2}{1-s}}\si^{(\frac{\tilde \si}{\mu_1-\rho_1})^\frac{C}{|\ln a|}}+C\k_0^2\tilde \si^2. \nonumber
\end{eqnarray}
By choosing $ \tilde \si=(\mu_1-\rho_1)|\ln \si(t_0)|^{-\frac{\ln a}{2C}} $, we have
\begin{equation}\label{14-14}
	\ba{ll}
	\ds \dbE \Vert u_i(t,\cd)\Vert_{L^\i(\cB_{\rho_{\overline k}}(w_k))}^2\\
	\ns\ds \le C{\k_0^{\frac{2}{1-s}}} |\ln \si |^{-\frac{\ln a}{C}}(\mu_1-\rho_1+e^{-|\ln \si|^{\frac{1}{2}}}|\ln \si|^\frac{\ln a}{C})\\
	\ns\ds \le C{\k_0^{\frac{2}{1-s}}}|\ln \si |^{-\frac{\ln a}{C}}.
	\ea 
\end{equation}
On the other hand, since $u_i\in C_\dbF([0,T]; L^2(\O; C^1(G_i(t)))$,  we have
\begin{equation}\label{14-15}
	\dbE |u_i(t_0,x_0)|^2\le \dbE \Vert u_i(t,\cd)\Vert_{L^\i(\cB_{\rho_{\overline k}}(w_k))}^2+{C\k_0^2}\rho_{\overline k}^2\le  \dbE \Vert u_i(t,\cd)\Vert_{L^\i(\cB_{\rho_{\overline k}}(w_k))}^2+C{\k_0^2}\tilde \si^2. 
\end{equation}
Combining (\ref{14-14}) and (\ref{14-15}), noting that $ \ln a \sim \frac{1}{2E}\eta_1\a, ~ (\a\to 0)$,  we yield (\ref{14-2}).\endpf

If we consider  the case of $ u=u_1-u_2 $ in $\cG((0,T))\= G_1((0,T))\cap G_2((0,T)) $,  we do not need to assume $ t_0-4\eta_1^{-1}\rho_1 \ge 0$ when we apply  Theorem \ref{mth}. So after removing the assumption $ \sqrt{t_0}-\frac{\rho_0\sin \a}{E(1+\sin \a)}\ge 0 $, we can still get (\ref{15-3}) and (\ref{14-9}). We introduce the following result and omit the proof.
\begin{theorem}\label{theorem4.3}
	Suppose Assumptions \ref{assumption2} and \ref{assumption1} hold.	Let $ t_0\in (0,T] $ and $ x_0\in \pa \cG(t_0) $ and suppose $ \pa\cG(t_0) $ is of Lipschitz class with the constants $ \rho_0$ and  $\a $.  Set
	\begin{equation}\label{23-5}
		\mu_1=\frac{\rho_0}{1+\sin\a},\q w_1=x_0+\mu_1\zeta,\q \cl \rho_1=\frac{1}{4E}\eta_1\mu_1\sin \a,\q \cl \si=\dbE \int_{\cB_{\cl\rho_1}(w_1)}u(t_0,x)^2dx.
	\end{equation}
	If $ \cl \si\in (0,1) $, then we have
	\begin{equation}\label{23-6}
		\dbE |u(t_0,x_0)|^2\le C\k_0^C|\ln \cl \si|^{-\frac{\a}{C}}.
	\end{equation}
\end{theorem}
\begin{remark}
	From Assumption \ref{assumption1}, we immediately yield $ \pa G_i(t) $ $ (i=1,2) $ is of Lipschitz class of $ R_0, \frac{\pi}{4} $. However, this may not holds for $ \cG(t)=G_1(t)\cap G_2(t) $.
\end{remark}

\section{Proofs of Theorems \ref{uniqueness} and  \ref{theorem1} }

Let $ u_1, u_2 $ be two solutions of equation (\ref{CD}) with two domains $ G_1((0,T)), ~G_2((0,T)) $ respectively and let  $I_1((0,T))$ and $I_2((0,T))$ be the unknown boundaries.  Denote $ \cG(\cd )\= G_1(\cd)\cap G_2(\cd) $.

\par {\bf Proof of Theorem \ref{uniqueness}. } Let $u=u_1-u_2$. Then $u$ solves
\begin{equation}\label{u1-u2}
	du-\D udt=(a_1\cd \n u+b_1 u)dt+c_1udW(t)\q \mbox{ in } \cG ((0,T)).
\end{equation}

Set
$$q_2\=\frac{1}{2}\Vert a_1\Vert_{L^{\infty}_{\dbF}(0,T;C^{n}(O;\dbR^n))}^2 +  \Vert b_1\Vert_{L^{\infty}_{\dbF}(0,T;C^{n}(O))} + \frac{1}{2}\Vert c_1\Vert_{L^{\infty}_{\dbF}(0,T;C^{n}(O))}^2+ 1.
$$ 
By Theorem \ref{SSUCP} and (\ref{o1}), we find $u\equiv 0$ in $G_1((0,t_0))\cap  G_2((0,t_0))$. Then we obtain that $u_1=u_1-u_2=0$ on $I_2 ((0,t_0))$.   Recalling  $I_1(0)=I_2(0)$,  multiplying the equation of $ u_1$ by $ e^{-q_2t} u_1 $ and integrating on $G_1((0,t_0))\setminus G_2((0,t_0))$, we have
\begin{eqnarray}\label{1}
	&& \dbE \int_{G_1((0,t_0))\setminus G_2((0,t_0))} e^{-q_2t}u_1(du_1-\D u_1dt)dx \nonumber\\	
	&& = \frac{1}{2}\dbE \int_{G_1(t_0)\setminus G_2(t_0)} e^{-q_2t_0} u_1^2(t_0, x)dx+q_0\dbE \int_{G_1((0,t_0))\setminus G_2((0,t_0))}e^{-q_2t_0}u_1^2dxdt\\
	&& \q +\dbE\int_{G_1((0,t_0))\setminus G_2((0,t_0))} \(-\frac{1}{2} e^{-q_0t}c_1^2u_1^2+e^{-q_0t}|\n u_1|^2\)dxdt \nonumber\\
	&& = \dbE\int_{G_1((0,t_0))\setminus G_2((0,t_0))} e^{-q_0t} u_1(a_1\cd\n u_1+b_1 u)dxdt. \nonumber
\end{eqnarray}
Then we find that
\begin{equation*}\label{19-3}
	\frac{1}{2}\dbE \int_{G_1(t_0)\setminus G_2(t_0)} e^{-q_2t_0} u_1^2(t_0, x)dx\le 0,
\end{equation*}
which implies $u_1(t_0, \cd )=0$ on $G_1(t_0)\setminus G_2(t_0)$. 

By a similar argument, we get $u_2(t_0, \cd )=0$ on $G_2(t_0)\setminus G_1(t_0)$.  

Combining with the SUCP  in Corollary \ref{SSUCP}, we obtain that $ u_i(t_0, \cd)=0 $ on $ G_i(t_0) $ for $ i=1,2 $,  provided $ G_1(t_0)\neq G_2(t_0) $. This contradicts Assumption \ref{assumption4} and completes the proof of Theorem \ref{uniqueness}.
\endpf 

\ms

\par Now   we are in the position to  prove Theorem \ref{theorem1}.  Borrowing some ideas from \cite{SV}, we need the following   propositions as preparation and we put the proofs in Appendix. 

\begin{proposition}\label{theoremC}
	Under Assumptions \ref{assumption2}--\ref{assumption3},  we have
	\begin{equation}\label{s1}
		\sup_{t_0\in[0,T]}\dbE \int_{G_i(t_0)\setminus \cG(t_0)} u_i(t_0,x)^2dx\le  \k_0^C \big(\ln |\ln \tilde \e|\big)^{-\frac{1}{n}} ,\qq i=1,2,\q \tilde \e\in (0,e^{-1}),
	\end{equation}
	where $ \tilde \e $ represent the observation error given in (\ref{e}).
\end{proposition}
\begin{proposition}\label{theoremC2}
	Suppose Assumptions \ref{assumption2}--\ref{assumption3} hold, and  $ \pa \cG(t) $ is of Lipschitz class with the uniform constants $ \rho_0\le \min\{R_0, \frac{1}{\sqrt{2e}}\}$ and  $\a $ for every $ t\in [0, T] $. There exist a constant $ C $ depending on $ \rho_0, \a  $ such that for any $ \tilde \e\in (0,1) $, we have
	\begin{equation}\label{14-16}
		\sup_{t_0\in[0,T]}\int_{G_i(t_0)\setminus \cG(t_0)} u_i(t_0,x)^2dx\le  \k_0^C|\ln \tilde \e|^{-\frac{1}{C}},\qq i=1,2,
	\end{equation}
	where $ \tilde \e $ represent the observation error given in (\ref{e}).
\end{proposition}

Let $ G\subset \dbR^n $, we introduce the following notation:
\begin{equation}\label{s4}
	(G)_{\d}\= \big\{x\in G\big|~ \cB_{\d}(x)\subset G \big\}.
\end{equation}

The following result is  also useful. 
\begin{proposition}\label{theoremC3}
	Suppose  Condition \ref{assumption4}   and Assumptions \ref{assumption2}-\ref{assumption3} hold. Let $ u_i $ be  the solution of (\ref{CD}) with the unknown boundary $ I_i(t) $  ~($i=1,2$). For every $ \overline \rho >0$,   $ t_0\in [0,T] $ and    $z\in (G_i(t))_{2\overline \rho}  $, we have
	\begin{equation}\label{16-1}
		\dbE \int_{\cB_{\overline \rho}(z)} u_i(t_0,x)^2dx\ge \frac{1}{C}\exp\(-C\ln \k_0 e^{\frac{C}{\min\{\overline \rho^n, t_0^{n/2}\}}}\). 
	\end{equation}
\end{proposition}

Next, we  introduce the definition of a modified distance.

Let $ G_1 $ and $ G_2 $ be bounded domains in $ \dbR^n $. The modified distance between $ G_1 $ and $ G_2 $ is given by
\begin{equation}\label{17-1}
	{\bf d}_m(G_1,G_2)=\max\Big\{\sup_{x\in \pa G_1} \dist (x, \cl G_2), \;\sup_{x\in \pa G_2}\dist (x, \cl G_1)\Big\}.
\end{equation}

We recall the following known result.
\begin{lemma}\label{pro1}\cite[Proposition 3.6]{GAEBERSV}
	Under Assumption \ref{assumption1}, there exist positive constants $\rho_0$, $\a$ and $d_0 $ such that if 
	\begin{equation}\label{17-2}
		{\bf d}(G_1(t), G_2(t))\le d_0,
	\end{equation}
	then there exists a  positive constant $ C $ such that
	\begin{equation}\label{17-3}
		{\bf d}(G_1(t), G_2(t))\le C {\bf d}_m(G_1(t), G_2(t))
	\end{equation}
	and $ \cG(t)=G_1(t)\cap G_2(t) $ is has boundary of Lipschitz class with constants $ \rho_0,\a  $.
\end{lemma}

\par {\bf Proof of Theorem \ref{theorem1}. } For the sake of simplicity we denote, for any $ t\in (0,T] $, $ {\bf d}(t)={\bf d}(G_1(t), G_2(t)) $ and $ {\bf d}_m(t)={\bf d}_m(G_1(t), G_2(t)) $. Set
\begin{equation}\label{17-4}
	\si=\max\Big\{\sup_{t\in [0,T]}\dbE \int_{G_1(t)\setminus \cG(t)} u_1(t,x)^2dx, \; \sup_{t\in [0,T]}\dbE  \int_{G_2(t)\setminus \cG(t)} u_2(t,x)^2dx\Big\}.
\end{equation}
The proof is divided into four steps. 

\ms

\par {\bf Step 1.}  In this step, we prove that 
\begin{equation}\label{17-5}
	{\bf d}_m(t_0)\le  \exp (\g(t_0)^C)\si^{\g(t_0)^{-C}},\q t_0\in(0,T].
\end{equation}

Let us assume, without loss of generality, that $ \si<1 $ and there exist $ x_0\in I_1(t_0)\subset \pa G_1(t_0) $ such that $ \dist(x_0, G_2(t_0))={\bf d}_m(t_0) $. By (\ref{17-4}),  we have
\begin{equation}\label{17-6}
	\int_{G_1(t)\cap \cB_{{\bf d}_m(t_0)}(x_0)} u_1(t_0, x)^2dx\le \si.
\end{equation}
Let   $ R(t_0)=\min\{ \sqrt t_0, R_0, \frac{1}{\sqrt{2e}}\} $. We deal with  two cases. 

\ss

Case (i).   ${\bf d}_m(t_0)\ge \eta_1 R(t) $. Apply Proposition \ref{theoremC3} by choosing $ \cl \rho=\frac{\eta_1R(t_0)}{3} $ and $ z=x_0-\frac{2\eta_1R(t_0)}{3}\nu  $, combining with (\ref{17-6}), we find that
\begin{equation}\label{17-7}
	\si\ge\dbE \int_{\cB_{\cl\rho}(z)} u_1(t_0,x)^2dx\ge \frac{1}{C}\exp\(-C\ln \k_0e^{\frac{C}{t_0^{n/2}}}\).
\end{equation}
Noting that ${\bf d}_m(t_0) $ is bounded, from (\ref{17-7}), we immediately obtain that
\begin{equation}\label{17-8}
	{\bf d}_m(t_0)\le C\exp\( C\ln \k_0 e^{\frac{C}{t_0^{n/2}}}\)\si. 
\end{equation} 

\ss

Case (ii).    ${\bf d}_m(t_0)< \eta_1R(t_0) $. Apply  Theorem \ref{main_theorem*} with $ r={\bf d}_m(t_0)$, $\rho=\eta_1R(t_0)$ and $R=R(t_0) $, we get that
\begin{equation}\label{17-9}
	\dbE \int_{G_1(t_0)\cap \cB_{\eta_1 R(t_0)}(x_0)}u_1(t_0,x)^2dx\le C\k_0^2 \si^{\frac{1}{C(\ln R(t_0)-\ln {\bf d}_m(t_0))}}.
\end{equation}
Meanwhile, similar to the argument in proving (\ref{17-7}), we find that
\begin{equation}\label{17-10}
	\dbE \int_{G_1(t_0)\cap \cB_{\eta_1 R(t_0)}(x_0)}u_1(t_0,x)^2dx\ge \frac{1}{C}\exp\(-C\ln \k_0 e^{\frac{C}{ t_0^{n/2}}}\).
\end{equation}
Combining (\ref{17-9}) and (\ref{17-10}), we get
\begin{equation}\label{17-11}
	\frac{1}{C(\ln R(t_0)-\ln {\bf d}_m(t_0))}\ln \si\ge -\ln C\k_0^2 -Ce^{\frac{C}{ t_0^{n/2}}}.
\end{equation}
Consequently,
\begin{equation}\label{24-7}
	{\bf d}_m(t_0)\le R(t_0)e^{\frac{\ln \si}{C}\[\ln C\k_0^2+Ce^{\frac{C}{t_0^{n/2}}}\]^{-1}}\le R(t_0) \si^{\[C(\ln \k_0)^2+e^{\frac{C}{t_0^{n/2}}}\]^{-1}}.
\end{equation}
From (\ref{17-8}) and (\ref{24-7}), we get (\ref{17-5}).

\ms

\par {\bf Step 2. } In this step, we  prove the following estimate:
\begin{equation}\label{18-1}
	{\bf d}(t_0)\le \exp (\g(t_0)^C)\si^{\g(t_0)^{-C}}.
\end{equation}

Without loss of generality, we may assume that there exists $ y_0\in \cl G_1(t_0) $ such that $ \dist (y_0, \cl G_2(t_0)) ={\bf d}(t_0)$. Let $ \d(t)=\dist(y_0, \pa G_1(t_0)) $. We cosider the following three cases:

\ss

Case (i) $ \d(t_0)\le \frac{1}{2}{\bf d}(t_0) $;

\ss

Case (ii) $ \frac{1}{2}{\bf d}(t_0)\le \d(t_0)\le \frac{1}{2}d_0 $, where $ d_0 $ given in Proposition \ref{pro1};

\ss

Case (iii) $ \d(t_0)> \max\{\frac{1}{2}{\bf d}(t_0), \frac{1}{2}d_0\} $.

\ss

\par   Case (i). Let $ z_0\in \pa G_1(t_0) $ be such that $ |y_0-z_0|=\d(t) $. We have
\begin{equation}\label{18-2}
	{\bf d}_m(t)\ge \dist (z_0, \cl G_2(t_0))\ge \dist (y_0, \cl G_2(t_0))-|y_0-z_0| = {\bf d}(t_0)-\d(t_0)\ge \frac{1}{2}{\bf d}(t_0).
\end{equation}
Therefore, by (\ref{17-5}), we get (\ref{18-1}).

Case (ii). By ${\bf d}(t_0)<d_0 $, Lemma \ref{pro1} and (\ref{17-3}), we have that ${\bf d}(t_0)\le C{\bf d}_m(t_0) $, which, together with (\ref{17-5}), implies (\ref{18-1}).

Case (iii). Let $ R_1(t_0)\=\min\{\frac{d_0}{2}, R(t_0)\} $ and $ {\bf d}_1(t_0)\=\min\{\frac{{\bf d}(t_0)}{2}, \eta_1R(t_0)\} $. By Assumption \ref{assumption1}, it is easy to check that
\begin{equation}\label{18-3}
	\cB_{{\bf d}_1(t_0)}(y_0)\subset G_1(t_0)\setminus \cl G_2(t_0)
\end{equation}
and
\begin{equation}\label{18-4}
	(t_0-R_1(t_0)^2,t_0]\t\cB_{R_1(t_0)}(y_0)\subset (0,t_0]\t G_1(t).
\end{equation}
By choosing $ r={\bf d}_1(t_0) $, $ \rho=\eta_1R_1(t_0) $ and $ R=R_1(t_0) $ in Theorem \ref{main_theorem*},   we get
\begin{equation}\label{18-5}
	\dbE \int_{\cB_{\eta_1R_1(t)}(y_0)} u_1(t_0,x)^2dx\le C\k_0^2\si^{\frac{1}{C(\ln R_1(t_0)-\ln {\bf d}_1(t_0))}}.
\end{equation}
On the other hand,  by choosing $ \cl \rho=\frac{\eta_1R_1(t_0)}{2} $ and $ z=y_0  $ in Proposition \ref{theoremC3}, we obtain 
\begin{equation}\label{18-6}
	\dbE  \int_{\cB_{\eta_1R_1(t_0)}(y_0)}u_1(t_0,x)^2dx\ge \frac{1}{C}\exp\(-C\ln \k_0e^{\frac{C}{ t_0^{n/2}}}\).
\end{equation}
If ${\bf d}_1(t_0)=\frac{{\bf d}(t_0)}{2} $, combining with (\ref{18-5}) and (\ref{18-6}), we get
\begin{equation}\label{18-7}
	{\bf d}(t_0)\le CR(t_0)\si^{\[C(\ln \k_0)^2+e^{\frac{C}{t_0^{n/2}}}\]^{-1}}.
\end{equation}
If ${\bf d}_1(t_0)=\eta_1R(t_0) $, by (\ref{18-6}), we see
\begin{equation}\label{18-8}
	\si \ge \frac{1}{C} \exp\(-C\ln \k_0 e^{\frac{C}{ t_0^{n/2}}}\),
\end{equation}
which yields
\begin{equation}\label{18-9}
	{\bf d}(t_0)\le C\exp \(C\ln \k_0 e^{\frac{C}{ t_0^{n/2}}}\)\si. 
\end{equation}
From (\ref{18-8}) and (\ref{18-9}), we obtain (\ref{18-1}).

\ms

{\bf Step 3. } In this step,  we will combine estimate (\ref{18-1}) with the stability estimate  given in Propositions \ref{theoremC} and \ref{theoremC2} to get the desired stability estimate.

\ss

From Proposition \ref{theoremC} and (\ref{18-1}), we find that
\begin{equation}\label{18-10}
	{\bf d}(t_0)\le \d_1 (t_0, \tilde \e)\= R(t_0)\exp (\g(t_0)^C)|\ln |\ln \tilde \e||^{-\g(t_0)^{-C}}.
\end{equation}
Notice that $ \d_1(t_0,\tilde \e) $ is increasing in $ \tilde \e\in (0,e^{-1}) $. There exist a constant $ \e_0 >0$ solving $ \d_1(t_0, \e_0)=d_0 $ such that  when $ \tilde \e<\e_0  $, we have ${\bf d}(t_0)\le \d_1(t_0, \tilde  \e)\le d_0 $. By Lemma \ref{pro1}, we have $ \pa \cG(t_0) $ is of Lipschitz class with the constants $ \rho_0$ and $\a $. By Proposition \ref{theoremC2}, we find that
\begin{equation}\label{18-11}
	{\bf d}(t_0)\le \exp (\g(t_0)^C)|\ln \tilde \e|^{-\g(t_0)^{-C}}.
\end{equation}
On the other hand, if $ \e_0\le \tilde \e< 1 $,  by $ \d_1(t_0, \e_0)=d_0 $, we see that
\begin{equation}\label{24-8}
	\ln |\ln \e_0|= d_0^{-\g(t_0)^C}\exp(\g(t_0)^{2C}).
\end{equation}
Consequently,
\begin{equation}\label{18-12}
	{\bf d}(t_0)\le C\(\frac{|\ln \e|}{|\ln \e_0|}\)^{-\g(t_0)^{-C}}\le \exp(e^{\g(t_0)^C})|\ln \tilde \e|^{-\g(t_0)^{-C}}.
\end{equation}
From (\ref{18-11}) and (\ref{18-12}), we get the desired result.
\endpf

\section{Appendix}

\par {\bf Proof of Proposition \ref{theoremC}. } We divide the proof into three steps.

\ss

{\bf Step 1}. Let $ \d <\min\{R_0, \frac{1}{\sqrt{2e}},\min_{t\in[0,T]}\dist(O_0,\pa G(t)) \} $, where $ R_0 $ is given in Assumption \ref{assumption1}. Set $ u=u_1-u_2 $. 
Since $ u(0, \cd) =0$,    for any $ t\in (0,T] $ and $ x_0\in (\cG(t))_{\frac{\d}{2}}  $,  we have $ (\max\{t_0-R^2, 0\}, t_0)\t \cB_{\frac{\d}{2E}}(x_0)\subset \cG((0,t_0)) $ (where $ E $ appears in Assumption \ref{assumption3}). Choosing $ R=\frac{\d}{2E} $, $ \rho=\frac{\eta_1\d}{2E}  $ and $ r=\frac{\eta_1\d}{6E} $  in Theorem \ref{mth}, we know there exists a  constant $ s \in (0,1)$, independent of $ x_0$ and $t $, such that
\begin{eqnarray}\label{s5}
	&& \dbE \int_{\cB_{\frac{\eta_1\d}{2E}}(x_0)}u(t,x)^2dx \nonumber\\
	&& \le C |\ln \d|^3\d^{n(1-s)} \k_0^{2-s} \(\dbE \int_{\cB_{\frac{\eta_1\d }{6E}}(x_0)} u(t,x)^2dx\)^{s}\\
	&& \le C\k_0^2\(\dbE \int_{\cB_{\frac{\eta_1\d }{6E}}(x_0)} u(t,x)^2dx\)^{s}. \nonumber
\end{eqnarray}
Let $ |x_1-x_0|=\frac{\eta_1\d}{3} $ and we immediately have $ \cB_{\frac{\eta_1\d}{6E}}(x_1)\subset \cB_{\frac{\eta_1\d}{2E}}(x_0) $.  Then we get
\begin{equation}\label{23-4}
	\dbE \int_{\cB_{\frac{\eta_1\d}{6E}}(x_1)}u(t,x)^2dx\le C\k_0^2\(\dbE \int_{\cB_{\frac{\eta_1\d }{6E}}(x_0)} u(t,x)^2dx\)^{s}.
\end{equation}
Irritating the inequality (\ref{s5})   $\lfloor\frac{C}{\d^n}\rfloor+1$ ($\lfloor\frac{C}{\d^n}\rfloor$ is the largest integer smaller than $\frac{C}{\d^n}$) times such that the left hand side of (\ref{23-4})  cover $ (\cG(t))_{\frac{\d}{2}} $,   assuming $ \d\le \dist(O_0, \pa G) $,  we find that
\begin{equation}\label{s6}
	\sup_{t\in [0,t_0]}\dbE \int_{(\cG(t))_{\frac{\d}{2}}} u(t,x)^2dx\le \frac{C\k_0^{\frac{2}{1-s}}}{\d^n}\tilde \e^{(\frac{1}{e})^{\frac{C}{\d^n}}}, \q t_0\in [0,T].
\end{equation} 

\ms

{\bf Step 2}. In this step, we prove the following inequality:
\begin{equation}\label{s6*}
	\dbE\int_0^{t_0} \Vert u(t,\cd)\Vert_{L^\i((\cG(t))_{\d})}^2 dt\le C
	\k_0^{\frac{2}{1-s}}\d^{-3n}\tilde \e^{(\frac{1}{e})^{\frac{C}{\d^n}}}, \q t\in [0,T].
\end{equation}

Choose a sequence of  cut-off functions $ \eta_i\in C^{\i}([0,T]\t \dbR^{n}) $ $ (i=1,2,\cdots, n) $ satisfying 
\begin{equation}\label{12-7}
	\eta_i(t,x)=\left\{
	\ba{ll}
	\ds 1,\q \mbox{in } (\cG(t))_{\frac{\d}{2}+\frac{i\d}{2n}},\\
	\ns\ds 0, \q \mbox{in } \cG(t)\setminus (\cG(t))_{\frac{\d}{2}+\frac{(i-1)\d}{2n}}.
	\ea 
	\right.
\end{equation}
Then, there exist a constant $ C>0 $ such that
$$\sum_{i=1}^n\big( |\eta_{it}|+\d |\n \eta_i|\big)\le C. $$
Multiply (\ref{CD}) by $ \eta_1^2 u $. Similar to the proof of the inequality \eqref{12-1},  we can prove that
\begin{equation}\label{12-8}
	\dbE \int_0^{t_0 }\int_{ (\cG(t))_{\frac{\d}{2}+\frac{\d}{2n}}}|\n u(t,x)|^2dxdt\le \frac{C}{\d^2} \dbE \int_0^{t_0 }\int_{ (\cG(t))_{\frac{\d}{2}}}u(t,x)^2dxdt.
\end{equation}
Next, we consider the following equation:
\begin{equation}\label{12-8*}
	du_{x_j}-\D u_{x_j}dt=\big(a_1\cd\n u+b_1u\big)_{x_j}dt + \big(c_{1}u \big)_{x_j}dW(t), \q j=1,2,\cdots, n.
\end{equation}
By multiplying (\ref{12-8*}) by $ \eta_2^2 u_{x_j} $, $j=1,\cds,n$, similar to (\ref{12-8}), using integration by parts, we  obtain 
\begin{equation}\label{h1}
	\ba{ll}
	\ds  \dbE\int_{0}^{t_0}\int_{\cG(t)_{\frac{\d}{2}+\frac{\d}{2n}}}\eta_2^2 u_{x_j}\big(du_{x_j}-\D u_{x_j}dt\big)dx \\
	\ns\ds =\dbE\int_{(\cG(t_0))_{\frac{\d}{2}+\frac{\d}{2n}}}\eta_2^2u_{x_j}^2(t_0,x)dx+ \dbE \int_{0}^{t_0}\int_{\cG(t)_{\frac{\d}{2}+\frac{\d}{2n}}} \big(-\eta_2\eta_{2t}u_{x_j}^2 +\eta_2^2 |\n u_{x_j}|^2\big)dxdt\\
	\ns\ds \q -\dbE \int_{0}^{t_0}\int_{\cG(t)_{\frac{\d}{2}+\frac{\d}{2n}}} \big[\eta_2^2  (c_{1x_j}u+c_1u_{x_j})^2+2\eta_2u_{x_j}\n \eta_2\cd \n u_{x_j}\big]dxdt.
	\ea 
\end{equation}
This, together with \eqref{12-8*}, implies that
\begin{equation}\label{h2}
	\dbE\int_{0}^{t_0}\int_{\cG(t)_{\frac{\d}{2}+\frac{\d}{2n}}}\eta_2^2|\n u_{x_j}|^2dxdt\le C\dbE\int_{0}^{t_0}\int_{\cG(t)_{\frac{\d}{2}+\frac{\d}{2n}}} \(u^2+\frac{1}{\d^2}|\n u|^2 \)dxdt.
\end{equation}

Combining (\ref{12-8}) and \eqref{h2},  we find that
\begin{equation}\label{12-9}
	\dbE \int_0^{t_0 }\int_{ (\cG(t))_{\frac{\d}{2}+\frac{\d}{n}}} |\n^2 u(t,x)|^2dxdt\le \frac{C}{\d^4} \dbE \int_0^{t_0 }\int_{ (\cG(t))_{\frac{\d}{2}}}u(t,x)^2dxdt.
\end{equation}
By repeating this process, we can finally get
\begin{equation}\label{12-10}
	\dbE \int_0^{t_0} \Vert u\Vert^2_{H^{\lfloor n/2 \rfloor +1}((\cG(t))_{\d})}dt\le \frac{C}{\d^{2n}}\dbE \int_0^{t_0 }\int_{ (\cG(t))_{\frac{\d}{2}}}u(t,x)^2dxdt.
\end{equation}
This, together with Sobolev embedding theorem, implies 
\begin{equation}\label{12-6}
	\dbE \int_0^{t_0}\Vert u(t,\cd)\Vert^2_{L^\i((\cG(t))_{\d})}dt\le \d^{-2n}\dbE \int_0^{t_0} \int_{(\cG(t))_{\d/2}} u(t,x)^2dxdt.
\end{equation}

Combining (\ref{12-6}) with (\ref{s6}), we conclude (\ref{s6*}).

\ss

{\bf Step 3}. In this step, we handle the term $  \dbE \int_{G_1(t_0)\setminus \cG(t_0)} u_1(t_0,x)^2 dx$. 

Notice that
\begin{equation}\label{s9}
	G_1(t_0)\setminus \cG(t_0)\subset \big[\big(G_1(t_0)\setminus (G_{1}(t_0))_{\d}\big)\setminus \cG(t_0)\big]\cup \big( (G_{1}(t_0))_{\d}\setminus (\cG(t_0))_{\d}\big).
\end{equation}
Recalling that $  u_1\in C_\dbF([0,T]; L^2(\O; C^1(G_1(t)))  $,  we have
\begin{equation}\label{s10}
	\dbE \int_{[G_1(t_0)\setminus (G_{1}(t_0))_{\d}]\setminus \cG(t_0)} u_1(t_0,x)^2dx\le C\k_0^2\d^{2+n},
\end{equation}
where we used Assumption \ref{assumption2} and  the fact that ${\bf d}\big(I_1(t_0), \big[G_1(t_0)\setminus (G_{1}(t_0))_{\d}\big]\setminus \cG(t_0)\big)\le C\d $.

Next, we denote $ (G_1\setminus \cG)_\d((0,t_0))\= \cup_{t\in(0,t_0)}\{t\}\t \big(\cG_1(t)\setminus \cG(t)\big)_{\d }$.  Multiplying the equation of $ u_1 $ by $ e^{-q_0t}u_1 $ ($ q_0 $ is given in (\ref{1})) and integrating on $ (G_1\setminus \cG)_\d((0,t_0)) $, we have
\begin{equation}\label{s11*}
	\ba{ll}
	\ds \dbE\int_{(G_1\setminus \cG)_\d((0,t_0))}e^{-q_2t}u_1 \big(du_1-\D u_1dt\big)dx\\
	\ns\ds =\frac{1}{2}\dbE\int_{(G_{1}(t_0))_{\d}\setminus (\cG(t_0))_{\d}} u_1(t_0,x)^2dx+q_0\dbE \int_{(G_1\setminus \cG)_\d((0,t_0))}e^{-q_2t_0}u_1^2dxdt\\
	\ns\ds\q  +\dbE\int_{(G_1\setminus \cG)_\d((0,t_0))} \(-\frac{1}{2} e^{-q_0t}c_1^2u_1^2 + e^{-q_0t}|\n u_1|^2\)dxdt \\
	\ns\ds  \q -\dbE\int_0^{t_0} \int_{\pa[(G_{1}(t))_{\d}\setminus (\cG(t))_{\d}]} u_1\n u_1\cd \nu_t(x)d\big[\pa\big[(G_{1}(t))_{\d}\setminus (\cG(t))_{\d}\big]\big]\\
	\ns\ds \q +\dbE\int_{\cup_{t\in(0,t_0)} \{t\}\t \pa (G_1(t)\setminus \cG(t))_\d } u_1^2\nu (t,x)\cd e_{1+n}d\big[\cup_{t\in(0,t_0)} \{t\}\t \pa \big(G_1(t)\setminus \cG(t)\big)_\d\big]\\
	\ns\ds  = \dbE\int_{G_1((0,t_0))\setminus G_2((0,t_0))} e^{-q_0t} u_1\big(a_1\cd\n u_1+b_1 u\big)dxdt,
	\ea 
\end{equation}
where $ \nu_t(x)$ (\resp $\nu(t,x) $) denotes the exterior unit normal to $ \pa[(G_{1}(t))_{\d}\setminus (\cG(t))_{\d}] $ (\resp $ \cup_{t\in(0,t_0)} \{t\}\t \pa (G_1(t)\setminus \cG(t))_\d $) and $ e_{1+n}=(1, 0,\cdots, 0) $. Then we get
\begin{equation}\label{s11}
	\ba{ll}
	\ds \dbE\int_{(G_{1}(t_0))_{\d}\setminus (\cG(t_0))_{\d}} u_1(t_0,x)^2dx\\
	\ns\ds  \le C\dbE\int_0^{t_0} \int_{\pa[(G_{1}(t))_{\d}\setminus (\cG(t))_{\d}]} u_1\n u_1\cd \nu_t(x)d\big[\pa\big[(G_{1}(t))_{\d}\setminus (\cG(t))_{\d}\big]\big]\\
	\ns\ds \q +C\dbE\int_{\cup_{t\in(0,t_0)} \{t\}\t \pa (G_1(t)\setminus \cG(t))_\d } u_1^2\nu (t,x)\cd e_{1+n}d\big[\cup_{t\in(0,t_0)} \{t\}\t \pa \big(G_1(t)\setminus \cG(t)\big)_\d\big].
	\ea 
\end{equation}
This, together with Assumption \ref{assumption2} and  H\"older's inequality, implies that
\begin{equation}\label{s12}
	\dbE \int_{(G_{1}(t_0))_{\d}\setminus (\cG(t_0))_{\d}} u_1(t_0,x)^2dx\le C\k_0\(\dbE \int_0^{t_0} \int_{\pa[(G_{1}(t))_{\d}\setminus (\cG(t))_{\d}]}u_1^2dtds\)^{1/2}.
\end{equation}
Clearly,
\begin{equation}\label{s13}
	\pa[(G_{1}(t))_{\d}\setminus (\cG(t))_{\d}]\subset \big\{\pa[(G_{1}(t))_{\d}\setminus (\cG(t))_{\d}]\cap \pa (G_1(t))_{\d}\big\}\cup \big\{\pa[(G_{1}(t))_{\d}\setminus (\cG(t))_{\d}]\cap \pa (G_2(t))_{\d}\big\}.
\end{equation}
Denote by $ \G_1(t) $ and $ \G_2(t) $ the part of $ \pa[(G_{1}(t))_{\d}\setminus (\cG(t))_{\d}]  $ in $ \pa (G_1(t))_{\d} $ and $ \pa (G_2(t))_{\d} $, respectively. For $ x\in \G_1(t) $, there exist $ y\in \pa G_1(t) $ such that
{\begin{equation}\label{s14}
		\dbE|u_1(t,x)|^2=\dbE |u_1(t,x)-u_1(t,y)|^2\le \dbE \Vert \n u_1\Vert_{L^\i(G_1(t))}^2|x-y|^2\le C\k_0^2 \d^2.
\end{equation}}
For $ x\in \G_2(t) $, recalling (\ref{s6*}), we have
\begin{equation}\label{s15}
	\ba{ll}
	\ds \dbE \int_0^{t_0} \int_{\G_2(t)} |u_1(t,x)|^2 d\G_2(t) dt \\
	\ns\ds \le 2 \dbE \int_0^{t_0}  \int_{\G_2(t)}\big( |u(t,x)|^2+ |u_2(t,x)|^2 \big)d\G_2(t) dt \\
	\ns\ds \le C\k_0^2 \d^2+ C\k_0^{\frac{2}{1-s}}\d^{-3n}\tilde \e^{(\frac{1}{e})^{\frac{C}{\d^n}}}.
	\ea 
\end{equation}
Set $ \d_1\=(2C)^{\frac{1}{n}} [\ln |\ln \tilde \e|]^{-\frac{1}{n}} $. If $  \d_1<\min \{R_0, \frac{1}{\sqrt{2e}}\} $, set $ \d=\d_1 $. By (\ref{s10}), (\ref{s13}) and (\ref{s15}),  we have
\begin{equation}\label{s16}
	\ba{ll}
	\ds \dbE 	\int_{G_1(t_0)\setminus \cG(t_0)} u_1(t_0,x)^2dx\\
	\ns\ds \le C\k_0\d+ C\k_0^{\frac{2-s}{1-s}}\d^{-\frac{3}{2}n}\tilde \e^{(\frac{1}{e})^{\frac{C}{\d^n}}}\\
	\ns\ds \le C\k_0^{\frac{2-s}{1-s}}[\ln |\ln \tilde \e|]^{-\frac{1}{n}}+C\k_0^{\frac{2-s}{1-s}}[\ln|\ln \tilde \e|]^{-\frac{1}{n}}[\ln|\ln \tilde \e|]^{\frac{3}{2}+\frac{1}{n}} e^{-\sqrt{|\ln \tilde \e|}}\\
	\ns\ds \le C\k_0^{\frac{2-s}{1-s}}[\ln |\ln \tilde \e|]^{-\frac{1}{n}}.
	\ea 
\end{equation} 
If $ \d_1\ge \min\{R_0, \frac{1}{\sqrt{2e}}\} $, then $ [\ln |\ln \tilde \e|]^{-\frac{1}{n}}\ge \frac{1}{C} $. Consequently,
\begin{equation}\label{s17}
	\dbE 	\int_{G_1(t_0)\setminus \cG(t_0)} u_1(t_0,x)^2dx\le C\k_0^2 [\ln |\ln \tilde \e|]^{-\frac{1}{n}}.
\end{equation}
Apply the same argument to $ u_2 $, we complete the proof.
\endpf

\ms

Before giving the detail for the proof of Proposition \ref{theoremC2}, we briefly illustrate the main ideas. Thanks to the extra regularity assumption on $ \pa \cG(t) $, compared with (\ref{s6}), we only need to irritate the two sphere one cylinder inequality by $\lfloor\frac{C}{\rho_0^n}\rfloor+1$ times. Then, by combing with the small propagation estimate in Theorem \ref{theorem4.3} and Assumption \ref{assumption2}, we can use similar technique in  the proof of Proposition \ref{theoremC} to get (\ref{14-16}).

\ms

\par {\bf Proof of Proposition \ref{theoremC2}. } By divergence theorem,  similar to the argument to prove  (\ref{s11}), we can show that,  for every $ t_0\in [0, T] $,
\begin{equation}\label{14-17}
	\ba{ll}
	\ds \dbE \int_{G_{1}(t_0)\setminus \cG(t_0)} u_1(t_0,x)^2dx\\
	\ns\ds \le C\dbE\int_0^{t_0} \int_{\pa[G_{1}(t)\setminus \cG(t)]} u_1\n u_1\cd \nu_t(x)d\big[\pa[G_{1}(t)\setminus \cG(t)]\big]\\
	\ns\ds +C\dbE\int_{\cup_{t\in(0,t_0)} \{t\}\t \pa (G_1(t)\setminus \cG(t)) } u_1^2\nu(t,x)\cd e_{1+n}d\big[\cup_{t\in(0,t_0)} \{t\}\t \pa (G_1(t)\setminus \cG(t))\big].
	\ea 
\end{equation}
Since for every $ t\in [0,T] $, 
\begin{equation}\label{14-18}
	\pa[G_{1}(t)\setminus \cG(t)] \subset  I_1(t)\cup  I_2(t),
\end{equation}
combing the fact that $ u_i(t,\cd)=0 $ on $  I_i(t) $, $ i=1,2 $, recalling $ u=u_1-u_2 $, we have
\begin{equation}\label{14-19}
	\ds \dbE \int_{G_{1}(t_0)\setminus \cG(t_0)} u_1(t_0,x)^2dx\le C\k_0\dbE \Vert u\Vert_{L^\i(\cup_{t\in(0,t_0)}\{t\}\t \{ I_2(t)\cap \pa[G_{1}(t)\setminus \cG(t)]\})}.
\end{equation}

Recalling \eqref{23-5} for $w_1$, for any $ x_0\in\pa \cG(t) $, we have  $ \cB_{\cl \rho_1}(w_1)\subset (\cG(t))_{\frac{3\rho_0}{8}} $. By applying Theorem \ref{theorem4.3} to $ u $ and noting (\ref{14-19}), we obtain that if $$ \sup_{t\in (0,T)}\linebreak \int_{(\cG(t))_{\frac{3\rho_0}{8}}}u(t,x)^2dx\in (0,1), $$
then
\begin{equation}\label{14-20}
	\dbE \int_{G_{1}(t_0)\setminus \cG(t_0)} u_1(t_0,x)^2dx\le C\k_0 \|\ln\(\sup_{t\in [0,t_0]} \int_{(\cG(t))_{\frac{3\rho_0}{8}}}u(t,x)^2dx\)\|^{-\frac{\a}{C}}.
\end{equation} 

\par On the other hand, recalling (\ref{s6}), by replacing $ \d $ by $ \frac{3}{4E}\rho_0 $ in (\ref{s5}), we  get there exist a absolute constant $ s\in (0,1) $, such that
\begin{equation}\label{15-1}
	\sup_{t\in [0,t_0]} \int_{(\cG(t))_{\frac{3\rho_0}{8}}}u(t,x)^2dx\le C\k_0^{\frac{2}{1-s}}\tilde \e^{\frac{1}{C}}, 
\end{equation}
where the generic constant $ C $ depends on $ \rho_0,\a  $.  Then, if $ C\k_0^{\frac{2}{1-s}}\tilde \e^{\frac{1}{C}}< \frac{1}{2} $,  we get
\begin{equation}\label{b11}
	\dbE \int_{G_{1}(t_0)\setminus \cG(t_0)} u_1(t_0,x)^2dx\le C|\ln \tilde \e|^{-\frac{1}{C}}.
\end{equation}
If $ \frac{1}{2}(C\k_0^{\frac{2}{1-s}})^{-C}\le \tilde \e<1$, we have 
\begin{equation}\label{b12}
	\dbE \int_{G_{1}(t_0)\setminus \cG(t_0)} u_1(t_0,x)^2dx\le\k_0^C|\ln \tilde \e|^{-\frac{1}{C}}.
\end{equation}
Combing  (\ref{b11}) and (\ref{b12}), with the same argument to $ u_2 $, we get the desired result. \endpf 

\ms

\par {\bf Proof of Proposition \ref{theoremC3}. } By Condition \ref{assumption4} and Assumption \ref{assumption1}, for any $ t_0\in [0,T] $, there exist $ x_0\in \G $ such that $ u_i(t_0,x_0)^2\ge F $.  From Assumption \ref{assumption2}, we know that $ \pa G_i(t) $ $ (i=1,2) $ is of Lipschitz class with the constants $ R_0$ and $\min\{\frac{\pi}{4}, \arcsin \frac{2\overline \rho}{3R_0}, \arcsin \frac{2E\sqrt t_0}{R_0}\} $.

Now we consider two cases. 

\ss

Case (i). $ \cl \rho\le \min\{R_0, 2E\sqrt{t_0}\} $. 

\ss

Noting that $ \frac{1}{4E}\eta_1 \frac{R_0\sin\a}{1+\sin \a}\le \frac{\eta_1}{4E} R_0 \frac{2\cl \rho}{3R_0}=\frac{\eta_1\cl \rho}{6E} $, by  applying Theorem \ref{theoremC1} on $ (t_0,x_0) $, we get
\begin{equation}\label{16-2}
	F\le C \k_0^C|\ln \si_i|^{-\frac{\a}{C}},
\end{equation}
where $ \a = \min\{\frac{\pi}{4}, \arcsin \frac{2\overline \rho}{3 R_0}, \arcsin \frac{2E\sqrt {t_0}}{R_0}\}$ and
\begin{equation}\label{16-3}
	\si_i\=\dbE \int_{\cB_{\frac{\eta_1\overline \rho}{6E}}(x_0+\frac{R_0\zeta}{1+\sin \a})} u_i(t_0,x)^2dx.
\end{equation}
Noting that $ \cB_{\frac{\cl \rho}{2}}(x_0+\frac{R_0\zeta}{1+\sin \a}) \subset G_i(t_0)$, by Assumption \ref{assumption3}, we have 
\begin{equation}\label{24-1}
	\(t_0-\frac{\cl\rho^2}{4E^2}, t_0\)\t \cB_{\frac{\cl\rho}{2E}}\(x_0+\frac{R_0\zeta}{1+\sin \a}\)\subset G_i( (0,t_0)).
\end{equation}
Applying the two sphere one cylinder inequality in  Theorem \ref{main_theorem*} with $ R=\frac{\cl \rho}{3E} $, $ \rho=\frac{\eta_1\cl \rho}{3E} $, $ r=\frac{\eta_1\cl \rho}{6E} $,  we obtain
\begin{equation}\label{24-2}
	\dbE \int_{\cB_{\frac{\eta_1\cl\rho}{3E}}(x_1)}u_i(t_0,x)^2dx\le C\k_0^{2-s}\(\dbE \int_{\cB_{\frac{\eta_1\cl\rho}{6E}}(x_1)}u_i(t_0,x)^2dx\)^s,
\end{equation}
where $ s=\frac{1}{C\ln 3\eta_1^{-1}} $ and $ x_1 $ be any point satisfying $ |x_1-x_0-\frac{R_0\zeta}{1+\sin \a}|=\frac{\eta_1\cl \rho}{6E} $. Noting the fact that
\begin{equation}\label{24-3}
	\cB_{\frac{\eta_1\cl\rho}{6E}}\(x_0+\frac{R_0\zeta}{1+\sin \a}\)\subset \cB_{\frac{\eta_1\cl\rho}{2E}}(x_1),
\end{equation} 
we have
\begin{equation}\label{24-4}
	\si_i \le C\k_0^{2-s}\(\dbE \int_{\cB_{\frac{\eta_1\cl\rho}{6E}}(x_1)}u_i(t_0,x)^2dx\)^s.
\end{equation}
Similar to the argument in (\ref{s5}), (\ref{23-4}), (\ref{s6}),  we can irritate (\ref{24-4}) by $\lfloor \frac{C}{\cl\rho^n}\rfloor+1 $ times to get
\begin{equation}\label{16-4}
	\si_i \le  C\k_0^{\frac{2}{1-s}} \(\int_{\cB_{\overline \rho}(z)}u_i(t_0,x)^2dx\)^{(\frac{1}{e})^{\frac{C}{ \cl\rho^n} }}.
\end{equation}
From (\ref{16-2}) and (\ref{16-4}), we see
\begin{equation*}\label{16-5}
	\frac{\k_0^{-C}}{C}e^{-(\frac{C}{F})^{\frac{C}{\a}}}\le \frac{\k_0^{-\frac{2}{1-s}}}{C}\si_i\le \(\int_{\cB_{\overline \rho}(z)}u_i(t_0,x)^2dx\)^{(\frac{1}{e})^{\frac{C}{ \cl\rho^n} }},
\end{equation*}
which implies 
\begin{equation}\label{24-5}
	\int_{\cB_{\overline \rho}(z)}u_i(t_0,x)^2dx\ge \exp\[ \(-C\ln \k_0-e^{(\ln C-\ln F)\frac{C}{\a}}\)e^{\frac{C}{\cl\rho^n}}\].
\end{equation}
Recalling that $ F $ is fixed, $\a = \min\{\frac{\pi}{4}, \arcsin \frac{2\overline \rho}{R_0}, \arcsin \frac{2E\sqrt t_0}{R_0}\}  $  and $ \cl\rho\le \min\{R_0, 2E\sqrt{t_0}\} $, we conclude that
\begin{equation}\label{16-6}
	\int_{\cB_{\overline \rho}(z)}u_i(t_0,x)^2dx\ge 	\frac{1}{C}\exp\(-C\ln \k_0e^{\frac{C}{\overline \rho^n}}\).
\end{equation}

On the other hand, if $ \cl \rho\ge\min\{ R_0, 2E\sqrt{t_0}\} $, by (\ref{16-6}), we find that
\begin{equation}\label{24-6}
	\int_{\cB_{\overline \rho}(z)}u_i(t_0,x)^2dx\ge \frac{1}{C}\exp\(-C\ln \k_0e^{\frac{C}{\min\{R_0^n, t_0^{n/2}\}}}\)\ge \frac{1}{C}\exp\[-C\ln \k_0e^{\frac{C}{\min\{\overline \rho^n, t_0^{n/2}\}}}\].
\end{equation}
Combining (\ref{16-6}) and (\ref{24-6}), we get (\ref{16-1}). 
\endpf


\begin{thebibliography}{99}
	
	\bibitem{GAEBERSV} G. Alessandrini, E. Beretta, E. Rosset and S. Vessella, \it Optimal stability for inverse elliptic boundary value problems with unknown boundaries. \sl  Ann. Scuola Norm. Sup. Pisa Cl. Sci., \rm  {\bf 29}(2000), 755--806.
	
	
	\bibitem{KBLC} K. Bryan and L. Caudill, \it Stability and reconstruction for an inverse problem for the heat equation. \sl Inverse Problems, \rm  {\bf 14}(1998), 1429--1453.
	
	\bibitem{BCERSV} B. Canuto, E. Rosset and S. Vessella, \it Quantitative estimates of unique continuation for parabolic equations and inverse initial boundary value problems with unknown boundaries, \sl Trans. Amer. Math. Soc., \rm  {\bf 354}(2002), 491--535.
	
	
	\bibitem{DCYOSP} D. Crisan, Y. Otobe and S. Peszat, \it Inverse problems for stochastic transport equations. \sl Inverse Problems, \rm  {\bf 31}(2015), Paper No. 015005.
	
	\bibitem{GDJZ} G. Da Prato and J. Zabczyk, \sl Stochastic equations in infinite dimensions.  \rm Cambridge University Press, Cambridge, 2014.
	
	\bibitem{FDWD} F. Dou and W. Du, \it Determination of solution of a stochastic parabolic equation by the terminal value. \sl Inverse Problems, \rm  {\bf 38}(2022), Paper No. 075010.
	
	\bibitem{JEGHMY} J. Elschner, G. Hu and M. Yamamoto, \it Single logarithmic conditional stability in determining unknown boundaries. \sl Appl. Anal., \rm   {\bf 99}(2020), 725--746.
	
	\bibitem{LEFJF} L. Escauriaza, F. J. Fern\'andez, \it Unique continuation for parabolic operators. \sl Ark. Mat., \rm 41(2003), 35--60.
	
	\bibitem{SVLE} L. Escauriaza and S. Vessella, \it Optimal three cylinder inequalities for solutions to parabolic equations with Lipschitz leading coefficients. Inverse problems: theory and applications. \sl Contemp. Math., \rm  {\bf 333}(2003), 79--87.
	
	
	\bibitem{Flandoli1990} F. Flandoli, \it Dirichlet boundary value problem for stochastic parabolic equations: compatibility
	relations and regularity of solutions. \sl Stochastics Stochastics Rep. \rm {\bf 29} (1990), 331--357.
	
	
	
	\bibitem{PK} P. Kotelenez, \it Stochastic ordinary and stochastic partial differential equations. \rm  Springer, New York, 2008.
	
	\bibitem{ZLQL} Z. Liao and Q. L\"u, \it Strong unique continuation property for stochastic parabolic equations. \rm https://doi.org/10.48550/arXiv.1701.02136. 
	
	\bibitem{JLLAMEZ} J. L\'imaco, L. A. Medeiros, E. Zuazua, \it Existence, uniqueness and controllability for parabolic equation in non-cylindrical domains. \sl Mat. Contemp.,  \rm {\bf 23}(2002), 49--70.
	
	
	
	
	\bibitem{QL1} Q. L\"u, \it Carleman estimate for stochastic parabolic equations and inverse stochastic parabolic problems, \sl Inverse Problems, \rm  {\bf 4}(2012), Paper No. 045008.
	
	\bibitem{Li2021} J. Li, H. Liu and S. Ma, \it Determining a random Schrödinger operator: both potential and source are random. \sl
	Comm. Math. Phys. \rm{\bf 381} (2021),  527--556.
	
	\bibitem{LY}  Q.~L\"u and Z.~Yin, \it  Unique continuation for
	stochastic heat equations. \sl ESAIM Control
	Optim. Calc. Var. \rm{\bf  21}(2015),  378--398.		
	
	\bibitem{LZ} Q. L\"u and X. Zhang, \it Global uniqueness for an inverse stochastic hyperbolic problem with three unknowns.
	\sl
	Comm. Pure Appl. Math.  \rm{\bf  68} (2015),  948--963.
	
	\bibitem{QLXZ} Q. L\"u and X. Zhang, \it Mathematical control theory for stochastic partial differential equations. \rm  Springer, Berlin, 2021.
	
	
	\bibitem{KDPGW} K. D. Phung and G. Wang, \it Quantitative unique continuation for the semilinear heat equation in a convex domain.  \sl J. Funct. Anal., \rm {\bf 259}(2010), 1230--1247.
	
	\bibitem{STXZ} S. Tang, X. Zhang, \it Null controllability for forward and backward stochastic parabolic equations, \sl SIAM J. Control Optim., \rm  {\bf 48}(2009), 2191--2216.
	
	
	\bibitem{SV} S. Vessella, \it Quantitative estimates of unique continuation for parabolic equations, determination of unknown time-varying boundaries and optimal stability estimates, \sl Inverse Problems, \rm  {\bf 24}(2008), Paper No. 023001.	
	
	\bibitem{MY} M. Yamamoto, \it Carleman estimates for parabolic equations and applications, \sl Inverse Problems, \rm  {\bf 25}(2009), Paper No. 123013, 75pp.
	
	\bibitem{Wu2020} B. Wu, Q. Chen and Z. Wang, Carleman estimates for a stochastic degenerate parabolic equation and applications to null controllability and an inverse random source problem. \sl
	Inverse Problems \rm {\bf 36} (2020),  075014, 38 pp.
	
	\bibitem{Wu2022} B. Wu, On the stability of recovering two sources and initial status in a stochastic hyperbolic-parabolic system.  \sl
	Inverse Problems \rm {\bf 38} (2022),   025010, 29 pp.
	
	\bibitem{GY} G. Yuan, \it Determination of two kinds of sources simultaneously for a stochastic wave equation. \sl
	Inverse Problems \rm {\bf 31} (2015),  085003, 13 pp.
	
	\bibitem{Yuan} G. Yuan, \it Conditional stability in determination of initial data for stochastic parabolic equations. \sl
	Inverse Problems \rm {\bf 33} (2017),  035014, 26 pp.
	
	
	
	\bibitem{XZ} X. Zhang, \it  Unique continuation for stochastic parabolic equations, \sl Differential Intergral Equations, \rm 21(2008), 81--93.
	
	
	
\end{thebibliography}
\end{document}